\documentclass[10pt]{amsart}
\usepackage{amsmath,amssymb,amsthm,amscd,mathrsfs}
\usepackage[all]{xy}
\usepackage{color} 
\CompileMatrices
\numberwithin{equation}{section}
\newtheorem{prop}{Proposition}[section]
\newtheorem{theo}[prop]{Theorem}
\newtheorem{lemm}[prop]{Lemma}
\newtheorem{coro}[prop]{Corollary}
\newtheorem{rema}[prop]{Remark}

\newtheorem{defi}[prop]{Definition}

\numberwithin{equation}{section}
\newcommand{\be}{\begin{equation}}
\newcommand{\ee}{\end{equation}}

\newcommand{\IP}{\mathbb{P}}%{{\relax{\rm I\kern-.18em P}}}

\newcommand\IZ{\mathbb {Z}}

\newcommand\IQ{\mathbb {Q}}
\newcommand{\IC}{\mathbb{C}}

\newcommand{\IR}{\mathbb{R}}

\newcommand{\ba}{\begin{array}}
\newcommand{\ea}{\end{array}}

\newcommand{\CX}{{\mathcal X}}

\newcommand{\IV}{{\mathbb V}}

\newcommand{\CS}{{\mathcal S}}

\newcommand{\IH}{{\mathbb H}} 
\newcommand{\bal}{\begin{aligned}}
\newcommand{\eal}{\end{aligned}}

\newcommand{\mmax}{{\mathrm{max}}}
\newcommand{\mmin}{{\mathrm{min}}}

\newcommand{\mfm}{{\mathfrak{M}}}

\newcommand{\longto}{\longrightarrow}
\newcommand{\lochom}{{\mathcal Hom}}

\newcommand{\CO}{{\mathcal O}}

\newcommand{\CE}{{\mathcal E}}

\newcommand{\CM}{{\mathcal M}}
\newcommand{\CL}{{\mathcal L}}

\newcommand{\CC}{{\mathcal C}}

\newcommand{\CQ}{{\mathcal Q}}

\newcommand{\dbar}{{\overline \partial}}
\newcommand{\obj}{{\mathfrak {Ob}}}

\title{Chamber Structure and Wallcrossing  in The ADHM Theory of Curves I}
\author{D.-E. Diaconescu}
\begin{document}

\begin{abstract} 
ADHM invariants are equivariant virtual invariants of moduli spaces 
of twisted cyclic representations of the ADHM quiver in the abelian category 
of coherent sheaves of a smooth complex projective curve $X$.  The goal of the present 
paper is to present a generalization of this construction employing 
a more general stability condition which depends on a real parameter. 
This yields a chamber structure in the ADHM theory of curves, residual ADHM 
invariants being defined by equivariant virtual integration in each chamber. 
Wallcrossing formulas will be presented in the second part of this work. 
\end{abstract} 
\maketitle 

\tableofcontents

\section{Introduction}\label{intro}

\subsection{Overview}\label{oversect}
Let $X$ be a smooth projective variety over $\IC$ equipped with a 
very ample line bundle $\CO_X(1)$. Let $M_1, M_2$ be fixed invertible 
sheaves 
on $X$ and let $E_\infty$ be a fixed coherent locally free 
$\CO_X$-modules. 
\begin{defi}\label{adhmsheaf}
An ADHM sheaf $\CE$ on $X$ with twisting data 
$(M_1,M_2)$ and framing data $E_\infty$ is a coherent $\CO_X$-module 
$E$ decorated by morphisms 
\[
\Phi_{i}:E\otimes_X M_i \to E, \qquad 
\phi:E\otimes_X M_1\otimes_X M_2 \to E_\infty, \qquad 
\psi:E_\infty \to E
\]
with $i=1,2$ satisfying the ADHM relation 
\be\label{eq:ADHMrelation}
\Phi_1\circ(\Phi_2\otimes 1_{M_1}) - \Phi_2\circ (\Phi_1\otimes 1_{M_2}) 
+ \psi\circ \phi =0.
\ee
An ADHM sheaf $\CE$ will be said to have Hilbert polynomial $P$ if 
$E$ has Hilbert polynomial $P$. If $X$ is a curve, an ADHM sheaf $\CE$ 
will be said to be of type $(r,e)\in \IZ_{\geq 0}\times \IZ$ if $E$ has 
rank $r\in \IZ_{\geq 0}$ and degree $e\in \IZ$.
\end{defi}

Motivated by string theoretic questions, the moduli problem for such decorated 
sheaves has been considered in detail in \cite{modADHM,GIT-decorated}. Very briefly, the main 
results obtained so far are as follows 
\begin{itemize}
\item[$\bullet$] According to \cite[Thm. 2.9.2.44, Thm. 2.9.2.5]{GIT-decorated} 
there exists a quasi-projective fine moduli space  of ADHM sheaves with fixed data 
$\CX=(X,M_1,M_2,E_\infty)$ and Hilbert polynomial $P$ subject to a certain 
stability condition.
The stability condition in question requires $E$ to be torsion free, 
$\psi$ to be nontrivial, and forbids the existence of nontrivial 
saturated proper subsheaves $0\subset E'\subset E$ 
preserved by $\Phi_1,\Phi_2$ and containing the image of $\psi$. 
This is essentially a cyclicity condition. 
\item[$\bullet$] If $X$ is a smooth projective curve over $\IC$, \cite[Thm. 1.2]{modADHM} 
shows that 
the above moduli space 
is equipped with a torus equivariant perfect tangent-obstruction theory
which yields a residual theory by virtual integration on the 
fixed loci. In particular the torus fixed loci have been shown to be proper
of finite type  over $\IC$. 
\item[$\bullet$] If $X$ is a smooth projective curve over $\IC$ and 
$E_\infty =\CO_X$, the resulting residual theory is equivalent to a 
local version of the stable pair theory constructed by Pandharipande and Thomas 
\cite{stabpairs-I}. The proof of this equivalence relies on the relative 
Beilinson spectral sequence \cite{orlov} for projective bundles. This was shown in 
\cite{modADHM} to yield a torus equivariant 
isomorphism of moduli spaces equipped
with perfect tangent obstruction theories between the moduli space 
of ADHM sheaves on $X$ 
and a moduli space of coherent systems on a projective 
bundle over $X$. The projective bundle $\pi:Y \to X$ is determined by the 
data $(M_1,M_2)$, that is $Y =\mathrm{Proj}(\CO_X \oplus M_1\oplus M_2)$. 
\end{itemize}

The main goal of the present paper is to 
construct a chamber structure in the ADHM theory of curves 
employing a more general stability condition which depends 
on a real stability parameter. Equivariant virtual invariants 
will be defined in each chamber and a wallcrossing result at the origin  
will be proven. 
Before summarizing the main results in more detail note that there at least 
two main reasons for this construction. 

First note that the presence of stability parameters 
is very natural in moduli problems for decorated sheaves 
\cite{verlinde,pairsI,pairsII,logflips,varmod,BDW,LePoitierA,LePoitierB,
stpairs,framed,Systcoh,equiv,augmented,triples,decorated, frHitchin,modquivers}. 
Recall that variations of the stability parameter have played a crucial 
role in the proof of the Verlinde formula \cite{verlinde,varmod}, birational 
transformations of moduli spaces \cite{verlinde,logflips,varmod}, 
the quantum cohomology of the grassmannian 
\cite{BDW}, as well as the topology of moduli spaces of sheaves 
of surfaces 
\cite{Flips,Systcoh,hodgenumbers,chalgsurf,chamber,GNYinst,NYblowupI,NYblowupII}.
In this paper it will be shown that the stability condition introduced in \cite{modADHM} 
for ADHM sheaves is in fact 
an asymptotic form of a more general $\delta$-stability condition, 
where $\delta\in \IR_{>0}$ 
 is a stability parameter. 
 
Moreover, moduli spaces of coherent systems have been 
recently related to curve counting problems by Pandharipande 
and Thomas \cite{stabpairs-I} and to generalized Donaldson-Thomas invariants 
by Joyce and Song \cite{genDTI}. Again, the stability condition employed in these 
cases is an asymptotic form of the 
stability conditions for coherent systems defined in \cite{LePoitierA,LePoitierB,Systcoh}. 
Therefore a natural question raised in the introduction of \cite{stabpairs-I}
is whether such virtual invariants admit 
deformations corresponding to variations of the stability 
parameter for coherent systems. This would require in principle the construction 
of a virtual cycle on the moduli space of stable coherent systems on smooth 
projective threefolds for arbitrary values of the stability parameter. As observed 
in the introduction of \cite{stabpairs-I}, such a 
cycle is not expected to exist in general, therefore the question seems to 
have a negative answer. However 
this paper will provide a positive 
answer for local stable pair theory 
theory, which was related to the ADHM theory of curves in \cite{modADHM}. 
It will be proven that such deformations of the local invariants exist, but the 
moduli spaces employed in the construction are not isomorphic to moduli 
spaces of coherent systems for generic values of the stability parameter
(see Remark (\ref{obsmod}) below). Note that a similar chamber structure in the stable pair theory of smooth projective 
Calabi-Yau threefolds has been previously constructed  in \cite{generating},
employing variations of stability conditions in the derived category\footnote{I thank A. Bayer 
for pointing this out.}. 

A different motivation for this paper resides in  string theoretic questions 
concerning wallcrossing phenomena for BPS states in Calabi-Yau 
compactifications and black holes \cite{DM-split,DM-crossing,GMN,JM}. More specifically, 
Jafferis and Moore \cite{JM} have shown that the spectrum of BPS 
states in a conifold compactification depends on an extra real 
parameter in addition to the expected complexified K\"ahler 
moduli. The constructions presented in this paper provide a rigorous mathematical framework 
for the wallcrossing considered in \cite{JM}. As shown in \cite{chamberII} they also yield 
a derivation of the halo wallcrossing formulas of Denef and Moore \cite{DM-split} 
for local rational curves. 

Finally note that similar constructions 
for Donaldson-Thomas invariants of toric resolutions of crepant threefold 
singularities have been carried out in \cite{szendroi-noncomm,MR-poisson,MR-noncomm,NH,
N-derived,wallcrossing} as well as the physics literature 
\cite{pyramid,OY}. 
These invariants are constructed in terms of moduli spaces of quiver 
representations, using the stability conditions defined in 
\cite{king}. Moreover, a generalization of the constructions employed in this paper has been
applied in \cite{stabquasimaps} to moduli spaces of stable quasimaps to holomorphic 
symplectic quotients.

The present work consists of two parts, the first being focused on existence and construction 
results as described below. The second part is concerned with wallcrossing formulas 
and some applications. 

\subsection{Construction results} 
This paper will consider 
ADHM sheaves on a smooth projective curve $X$ over $\IC$
with fixed twisting data $(M_1,M_2)$ and fixed framing data 
$E_\infty =\CO_X$.
Moreover, $M_1,M_2$ will be chosen so that 
there is an isomorphism $M_1\otimes_X M_2\simeq K_X^{-1}$.
This condition 
will be needed in the construction of a symmetric perfect obstruction 
theory for generic values of the stability parameter. 

The $\delta$-stability condition for ADHM sheaves is introduced 
in section (\ref{defbasic}), definition (\ref{deltastable}). Note that 
 lemma (\ref{duallemma}) allows us to restrict our treatment 
to positive values of the stability parameter. 
Several basic properties as well as preliminary boundedness results 
are proven in section (\ref{modspace}). 

Sections (\ref{absubcat}), (\ref{stabsubcat}) reformulate
 $\delta$-stability for ADHM sheaves 
as a slope stability condition in an abelian subcategory of quiver 
sheaves. This condition belongs to the class of slope stability conditions for quiver sheaves 
introduced by Alvarez-Consul and Garcia-Prada in 
\cite{dimred,HKquivers}. Similar stability conditions have been also considered in 
\cite{gauge-virtual}.
These results will be needed in the analysis of wallcrossing 
behavior in later sections. Some homological algebra results 
 are proven in section  (\ref{extsect}). 

Section (\ref{chambersect}) consists of a detailed analysis of variations of 
the stability parameter $\delta\in \IR_{>0}$. It is proven in section 
(\ref{crtstabpar}), lemmas (\ref{chamber}), 
(\ref{largedelta}), (\ref{finitecrt}),
that for fixed rank $r\in \IZ_{\geq 1}$ and fixed degree 
$e\in \IZ$ there are finitely many critical values $\delta_i\in \IR_{>0}$, 
$i=1,\ldots, N$ dividing the positive real axis into stability chambers. 
The set of $\delta$-stable ADHM sheaves of type $(r,e)$ is constant within 
each chamber, and strictly semistable objects may exist only if $\delta$ 
takes a critical value. Moreover, for sufficiently large $\delta$, $\delta$-stability 
is equivalent to the stability condition used in \cite{modADHM}. 
Section (\ref{wallsect}) is focused on the wallcrossing behavior of 
$\delta$-stable objects when $\delta$ specializes to a critical value. 
In particular lemmas (\ref{specstab}), (\ref{minusdestabobj}), respectively (\ref{plusminuslemm}), (\ref{zeroHN}) examine the interaction 
between generic $\delta$-stability and semistability at a critical 
value, respectively the origin. 

The moduli problem for $\delta$-semistable ADHM sheaves is the main subject of  
section (\ref{adhmmod}). Standard arguments imply that flat families of 
$\delta$-semistable ADHM sheaves form a groupoid $\mfm_\delta^{ss}(\CX,r,e)$ over the category of 
schemes over $\IC$. The first result proven in this section is the following.
\begin{theo}\label{modthm} 
The groupoid $\mfm_\delta^{ss}(\CX,r,e)$, where $\CX$ denotes the triple 
$(X,M_1,M_2)$,  is an Artin stack of finite 
type over $\IC$ for any $(r,e)\in \IZ_{\geq 1}\times \IZ$ and 
any $\delta \in \IR_{>0}$. If $\delta \in \IR_{>0}$ is noncritical
of type $(r,e)$, the moduli stack 
$\mfm_\delta^{ss}(\CX,r,e)$ is isomorphic to a quasi-projective scheme 
over $\IC$. 
\end{theo}

Other moduli stacks needed in the second part of this paper  are similarly constructed in section (\ref{biggerstack}). 

\begin{rema}\label{obsmod} 
According to lemma (\ref{largedelta}), the moduli space of 
$\delta$-semistable ADHM sheaves is isomorphic to the moduli space
constructed in \cite[Thm 1.1]{modADHM} for sufficiently large $\delta$. 
In particular it is also isomorphic to a quasi-projective moduli space
of asymptotically stable coherent systems as proven in 
\cite[Thm. 1.11]{modADHM}. The proof of this theorem relies on a 
vanishing result \cite[Lemma 2.5]{modADHM} 
for the morphism $\phi: E\otimes_X M_1\otimes_X M_2  \to \CO_X$ 
for asymptotically stable ADHM sheaves $\CE=(E,\Phi_1,\Phi_2,\phi,\psi)$. 
However, such a vanishing result no longer necessarily holds when $\delta$ 
lies in other chambers on the positive real axis. Therefore in this case, 
the moduli space of $\delta$-semistable ADHM sheaves is not expected to be 
isomorphic to a moduli space of $\delta$-semistable coherent systems.
\end{rema}

Next note that there is a 
natural torus ${\bf T}=\IC^\times\times \IC^\times$ action on the moduli 
space of $\delta$-stable ADHM sheaves, defined by scaling the ADHM 
data as follows
\[
(t_1,t_2)\times (E,\Phi_1,\Phi_2,\phi,\psi) \to (E, t_1\Phi_1,t_2\phi_2, t_1t_2
\phi, \psi). 
\]
Let ${\bf S}\simeq \IC^\times \subset {\bf T}$ denote the antidiagonal torus defined by the 
embedding $t\to (t^{-1},t)$. 
The next result proven in section (\ref{torvirt}) establishes the 
existence of a torus equivariant perfect tangent-obstruction theory
for noncritical values of the stability parameter. 
\begin{theo}\label{virsmooth}
Let $\delta\in \IR_{>0}$ be a noncritical stability parameter of type 
$(r,e)\in \IZ_{\geq 1}\times \IZ$. Then 
the  moduli scheme $\mfm_\delta^{ss}(\CX,r,e)$ has a 
{\bf T}-equivariant as well as ${\bf S}$-equivariant  
perfect tangent-obstruction theory. 
Moreover, 
the perfect tangent-obstruction theory of $\mfm_\delta^{ss}(\CX,r,e)$ 
is ${\bf S}$-equivariant symmetric. 
\end{theo}  

Using this structure, one would like to define a residual $\delta$-ADHM 
theory of curves by virtual integration on the torus fixed loci. In particular, 
a properness result for the fixed loci is needed. Properness of the 
{\bf T}-fixed loci can be proven by analogy with \cite[Prop. 3.15]{modADHM}. 
However, properness of the {\bf S}-fixed loci is more difficult,  and will require 
an inductive proof, which is presented in section (\ref{toractsect}). 
The final 
result is: 
\begin{theo}\label{Sfixed}
Let $\delta\in \IR_{>0}$ be a noncritical stability parameter of type 
$(r,e)\in \IZ_{\geq 1}\times \IZ$. Then the fixed locus 
$\mfm_\delta^{ss}(\CX,r,e)^{\bf S}$ is a projective 
scheme over $\IC$.
\end{theo}

An immediate corollary of theorem (\ref{Sfixed}) is 
\begin{coro}\label{Tfixed}
Under the conditions of theorem (\ref{Sfixed}), 
the fixed locus $\mfm_\delta^{ss}(\CX,r,e)^{\bf T}$ is also a projective scheme over $\IC$.
\end{coro} 

Theorems (\ref{virsmooth}) and (\ref{Sfixed}) and corollary (\ref{Tfixed}) 
imply via \cite{GP} the following 
\begin{coro}\label{indvircycle}
Suppose the conditions of theorem (\ref{Sfixed}) are satisfied. 
Then the following statements hold. 

$(i)$ The {\bf T}-fixed locus $\mfm_\delta^{ss}(\CX,r,e)^{\bf T}$ 
is equipped with an induced 
perfect tangent-obstruction theory, which yields a 
virtual fundamental cycle 
$$[\mfm_\delta^{ss}(\CX,r,e)^{\bf T}]\in 
A_{\bullet}(\mfm_\delta^{ss}(\CX,r,e)^{\bf T})$$
and a virtual 
normal bundle $N_{\mfm_\delta^s(\CX,r,e)^{\bf T}/ 
\mfm_\delta^{ss}(\CX,r,e)}\in K^0_{\bf T}( \mfm_\delta^{ss}(\CX,r,e)^{\bf T})$
in the equivariant $K$-theory of locally free sheaves of the fixed locus.

$(ii)$ Analogous results hold for the {\bf S}-fixed loci. Moreover, in this case, 
the induced perfect tangent-obstruction theory is symmetric and the 
resulting virtual fundamental cycle is a $0$-cycle. 
\end{coro}

Using corollary (\ref{indvircycle}) the residual ADHM theory 
of the data $\CX=(X,M_1,M_2)$ is defined as follows 
\begin{defi}\label{deltainv} 
Let $(r,e)\in \IZ_{\geq 1}\times \IZ$, and $\delta \in
\IR_{>0}$ be a noncritical stability parameter. Then 
the {\bf T}-equivariant 
$\delta$-ADHM invariant of type $(r,e)$ is defined by 
\[
A^{\bf T}_{\delta}(r,e) = \int_{[\mfm_\delta^{ss}(\CX,r,e)^{\bf T}]}
e_{\bf T}^{-1}(N_{\mfm_\delta^{ss}(\CX,r,e)^{\bf T}/ 
\mfm_\delta^{ss}(\CX,r,e)}).
\]
The {\bf S}-equivariant $\delta$-ADHM invariant of type $(r,e)$ 
$A^{\bf S}_{\delta}(r,e)$ is defined analogously. 
\end{defi}

Finally, the goal of section (\ref{versCS})  is to prove that 
the deformation results obtained by Joyce and Song for coherent sheaves on Calabi-Yau 
threefolds also hold  for locally free  ADHM sheaves on curves. 
More precisely Theorems (\ref{theoremA}), (\ref{theoremB}) 
are entirely analogous to Theorems \cite[Thm. 5.2]{genDTI}, 
\cite[Thm. 5.3]{genDTI} and follow by the same type of deformation 
theory arguments applied to locally free ADHM sheaves. 
Moreover it is proven that the Behrend function identities 
proven in  \cite[Thm. 5.9]{genDTI} also hold with appropriate 
modifications in the present case. 
Using these statements, the main results in the theory of generalized 
Donaldson-Thomas invariants of \cite{genDTI} apply to ADHM sheaf 
invariants, allowing one to derive explicit wallcrossing formulas. 
This will be presented in detail in the second part of this paper.  

{\it Acknowledgements}. 
I am very grateful to Arend Bayer,
Ugo Bruzzo, Wu-Yen Chuang, Daniel Jafferis, Jan Manschot, Greg Moore, 
Kentaro Nagao, Alexander Schmitt, Andrei Teleman, Yokinobu Toda, and especially 
Ron Donagi, Liviu Nicolaescu and Tony Pantev for very helpful discussions and 
correspondence. 
I owe special thanks to Tony Pantev for his patience, and to 
Alexander Schmitt for making a preliminary version 
of \cite{GIT-decorated} available to me prior to publication. 
I would also like to thank 
Ionut Ciocan-Fontanine, Bumsig Kim and Davesh Maulik for collaboration 
on a related project.  The partial support of NSF grants PHY-0555374-2006 
and NSF-PHY-0854757 is also 
acknowledged. 

{\bf Notation and Conventions}. 
Throughout this paper, we will denote by ${\mathfrak S}$ 
the category of schemes of finite type over 
$\IC$. For any such schemes $X,S$  
we  set $X_S = S\times X$ and $X_s=\mathrm{Spec}(k(s))\times_S X$ 
for any point 
$s\in S$, where $k(s)$ is the residual field of $s$. Let also 
$\pi_S: X_S \to S$, $\pi_X:X_S \to X$ denote the canonical projections. 
We will also set $F_S = \pi_X^*F$ for 
any $\CO_X$-module $F$. Given a morphism $f:S'\to S$, 
we will denote by $f_X = f\times 1_X:X_{S'}\to X_S$. 
Any 
morphism $f:S'\to S$, yields a  commutative diagram of the form 
\[
\xymatrix{
X_{S'} \ar[r]^-{p_X'} \ar[d]_-{f_X} & X \ar[d]^-{1_X} \\
X_S \ar[r]^-{\pi_X} & X. \\
}
\]
Then for any $\CO_X$-module $F$ there is a  
canonical isomorphism 
$F_{S'}\simeq f_X^* F_S$
which will be implicit in the following.

\section{$\delta$-stability for ADHM sheaves}\label{adhmstab}

\subsection{Definition and basic properties}\label{defbasic}
Let $X$ be a smooth projective curve over a $\IC$-field $K$
 equipped with
a very ample line bundle $\CO_X(1)$.
Let $M_1,M_2$ be fixed line bundles on $X$ so that  
$M_1\otimes_X M_2\simeq K_X^{-1}$. Such an isomorphism will 
be fixed throughout this paper, and we will also set $M=M_1\otimes_X M_2$. Note that the field of definition 
of $X$ is taken to be an extension of $\IC$ since such curves occur 
naturally as fibers in flat families over non-closed points. Therefore, 
in order 
to formulate the moduli problem, 
one has to define the basic objects of study over arbitrary extensions 
of $\IC$, by analogy with moduli of sheaves \cite{huylehn}.

As stated in the beginning of section (\ref{oversect}), we
 will consider ADHM sheaves \cite[Def. 2.1]{modADHM}
$\CE=(E,\Phi_1,\Phi_2,\phi,\psi)$ on $X$ 
with twisting data $(M_1,M_2)$ and framing data $E_\infty=\CO_X$.
Given a coherent locally free $\CO_X$-module 
$E$  we will denote by $r(E)$, $d(E)$, 
$\mu(E)$ the rank, degree, respectively slope of $E$.
An ADHM sheaf $\CE$ will be called locally free if $E$ is a coherent 
locally free $\CO_X$-module. 
If $\CE$ is a locally free ADHM sheaf 
the pair $(r(E), d(E))$ will be called the type 
of $\CE$.

Let $\delta\in \IR\setminus\{0\}$ be a stability parameter. For a nontrivial locally free ADHM 
sheaf  $\CE$ 
we define the $\delta$-slope of $\CE$ to be 
\[
\mu_\delta(\CE) = \mu(E) + {\delta \over r(E)}.
\]
Moreover a subsheaf $0\subset E'\subset E$ 
will be called $\Phi$-invariant if $\Phi_i( E'\otimes M_i)\subseteq 
E'$ for $i=1,2$. 

\begin{defi}\label{deltastable}
Let $\delta\in \IR\setminus\{0\}$ be a stability parameter. 
A nontrivial locally free ADHM sheaf $\CE=(E,\Phi_1,\Phi_2,\phi,\psi)$ is 
$\delta$-(semi)stable if the following conditions 
are satisfied 
\begin{itemize}
\item[$(i)$] $\psi$ is nontrivial if $\delta>0$ respectively  
$\phi$ is nontrivial if $\delta <0$.
\item[$(ii)$] 
Any $\Phi$-invariant nontrivial proper saturated 
subsheaf $0\subset E'\subset E$ 
so that $\mathrm{Im}(\psi) \subseteq E'$ satisfies 
\be\label{eq:slopecondA}
\mu(E')+ {\delta  \over r(E')}  \
(\leq)  \ \mu_\delta(\CE).
\ee
\item[$(iii)$] 
Any $\Phi$-invariant nontrivial proper saturated 
subsheaf $0\subset E'\subset E$ 
so that $E'\otimes_X M \subseteq \mathrm{Ker}(\phi)$ satisfies 
\be\label{eq:slopecondB} 
\mu(E') \
(\leq)  \ \mu_\delta(\CE).
\ee
\end{itemize}
 \end{defi}
 
We also define (semi)stability at $\delta=0$ as follows. 
\begin{defi}\label{zerostab} 
A nontrivial locally free  ADHM sheaf $\CE= (E, \Phi_{1,2}, \phi, \psi)$ be 
on $X$ is $0$-(semi)stable if any $\Phi$-invariant proper nontrivial saturated 
subsheaf 
$0\subset E'\subset E$ so that $\mathrm{Im}(\psi) \subseteq E'$ 
or $E'\otimes_X M \subseteq \mathrm{Ker}(\phi)$ satisfies 
\be\label{eq:zerostabcond}
\mu(E') \ (\leq ) \ \mu(E). 
\ee
\end{defi}
 
  Let $\CE=(E, \Phi_1,\Phi_2,\phi,\psi)$ be a locally free ADHM sheaf on $X$
of type $(r,e)\in \IZ_{\geq 1}\times \IZ$. 
Then the data 
\be\label{eq:dualADHM} 
\begin{aligned}
{\widetilde E} & = E^\vee \otimes_X M^{-1}\\
{\widetilde \Phi}_i & = (\Phi_i^\vee\otimes 1_{M_i}) \otimes 1_{M^{-1}}: {\widetilde E}
\otimes 
M_i \to {\widetilde E} \\
{\widetilde \phi} & = \psi^\vee \otimes 1_{M^{-1}} : {\widetilde E}\otimes_X {M} \to \CO_X
\\
{\widetilde \psi} & = \phi^\vee : \CO_X \to {\widetilde E} \\
\end{aligned}
\ee
with $i=1,2$, determines a locally free ADHM sheaf ${\widetilde \CE}$ of type 
$${\widetilde r} = r\qquad {\widetilde e} = -e +2r(g-1)$$ where $g\in \IZ_{\geq 0}$ is the 
genus of $X$. ${\widetilde \CE}$ will be called the dual of $\CE$ in the following. 
Note that ${\widetilde \CE}$ has the same twisting and framing data as $\CE$. 

\begin{lemm}\label{duallemma}
Let $\delta\in \IR_{>0}$ be a positive stability parameter and let 
$\CE$ be a locally free ADHM sheaf on $X$. 
Then $\CE$ is $\delta$-(semi)stable if and only if ${\widetilde \CE}$ is 
$(-\delta)$-(semi)stable.
\end{lemm}

{\it Proof.} Straightforward verification of stability conditions. 

\hfill $\Box$

Using lemma (\ref{duallemma}) it suffices to consider only positive 
stability parameters $\delta\in \IR_{>0}$ from this point on.

 \subsection{Boundedness results}\label{modspace}
This subsection consists of several boundedness results for semistable 
ADHM sheaves required at later stages in the paper. 
As in the previous subsection, 
$X$ is a smooth projective curve over a field $K$ over 
$\IC$ and $M_1, M_2$ are fixed line bundles on $X$ equipped with an isomorphism 
$M_1\otimes_X M_2\simeq K_X^{-1}$. 
\begin{lemm}\label{boundedness} 
Let 
$(r,e)\in \IZ_{\geq 1}\times \IZ$ be a fixed type. 
Then the set of isomorphism classes of locally free sheaves $E$ of type 
$(r,e)$ 
on $X$ so that 
$(E,\Phi_{1,2},\phi,\psi)$ is a $\delta$-semistable ADHM sheaf 
for some $\delta\in \IR_{\geq 0}$ and some 
morphisms $(\Phi_{1,2},\phi,\psi)$ is bounded. 
\end{lemm} 

{\it Proof.} 
The proof will be based on Maruyama's theorem \cite{Maruyama-bound}
\begin{theo}\label{Maruyama} 
$\mathrm{(Maruyama)}$. 
A family of torsion free sheaves $E$ with fixed Hilbert polynomial 
$P$ and $\mu_{\mmax}(E) \leq C$ for a fixed constant 
$C$ is bounded. 
\end{theo}

\noindent
and the following standard technical lemma 
(used for example in the proof of 
\cite[Prop. 3.2]{semistpairs}, \cite[Thm. 3.1]{projective}). 
 \begin{lemm}\label{HNfiltr} 
Let $E$ be a torsion-free sheaf of rank $r\geq 2$ on $X$. 
Suppose $E$ is not semistable, and let 
\[
0=HN_{0}(E)\subset HN_1(E) \subset \ldots \subset HN_h(E)=E
\]
be the slope Harder-Narasimhan filtration of $E$. 
Then 
\be\label{eq:HNineq}
\mu(HN_1(E)) + (r-1) \mu(E/HN_{h-1}(E)) \leq r \mu(E).
\ee
\end{lemm}
According to Maruyama's theorem it suffices to prove that 
there exists a constant $C$ independent on $\delta$ 
so that $\mu_{\textrm{max}}(E)\leq C$ 
for all $\delta$-semistable ADHM sheaves 
$\CE=(E,\Phi_{1,2},\phi,\psi)$ on $X$ 
of type $(r,e)$ and all $\delta\in \IR_{\geq 0}$.

If $E$ is semistable 
$\mu_{\textrm{max}}(E) =\mu(E)$ is clearly bounded.
In particular this is the case if $r=1$, hence 
we will assume $r\geq 2$ from now on in this proof. 
 Suppose $E$ is not semistable of rank $r\geq 2$, and let 
\be\label{eq:HNa} 
0= HN_0(E) \subset HN_1(E) \subset \cdots \subset HN_h(E) =E, 
\ee 
$h\geq 2$, be the Harder-Narasimhan filtration of $E$. 
Note that the successive quotients are locally free and semistable.  
In particular this implies that $h\leq r$.

Suppose first $\CE$ is $\delta$-semistable for some $\delta\in \IR_{>0}$
or that $\CE$ is 0-semistable and $\psi:E_\infty \to E$ is not identically zero. 
 Let  $j_\psi\in \{1,\ldots,h-1\}$ be the index determined by 
\[
\begin{aligned}
& \mathrm{Im}(\psi)\nsubseteq HN_{j}(E), \qquad \mathrm{for}\ j \leq j_\psi \\
& \mathrm{Im}(\psi) \subseteq HN_j(E),\qquad \mathrm{for}\ j \geq j_\psi +1.
\end{aligned}
\]
Note that the morphism 
\[
{\overline \psi}: E_\infty \to E/HN_{j_\psi}(E) 
\]
is nontrivial and \cite[Lemma 1.3.3]{huylehn} 
implies that
\[
\mu_{\mmin}(E_\infty) \leq \mu_{\mmax}(E/HN_{j_\psi}(E)). 
\]
By construction (see the proof of \cite[Thm. 1.3.4]{huylehn})
we have $\mu_{\mmax}(E/HN_{j_\psi}(E))= 
\mu(HN_{j_\psi+1}(E)/HN_{j_\psi}(E))$, therefore we obtain 
\be\label{eq:boundA}
\mu_{\mmin}(E_\infty) \leq \mu(HN_{j_\psi+1}(E)/HN_{j_\psi}(E)).
\ee
Moreover, if $j_\psi = h-1$, inequality \eqref{eq:boundA} specializes to 
\[
\mu_{\mmin}(E_\infty) \leq \mu(E/HN_{h-1}(E)).
\]
which yields 
\be\label{eq:boundB} 
-\mu(E/HN_{h-1}(E)) \leq - \mu_{\mmin}(E_\infty).
\ee
If $j_\psi < h-1$, we claim $HN_j(E)$ cannot be $\Phi$-invariant, 
for any $j_\psi\leq j \leq n-1$. 
According to the $\delta$-stability condition, if $HN_j(E)$ is
$\Phi$-invariant, $i=1,2$ for some  $j_\psi\leq j \leq h-1$,
it follows that 
\[
\mu(HN_j(E)) +{\delta \over r(HN_j(E))}\leq  
\mu(E) + {\delta \over r}
\]
Since $\delta\geq 0$ and $r(HN_j(E))<r$, this yields a contradiction because 
$\mu(HN_j(E))>\mu(E)$ for all $j=1,\ldots, h-1$.
This proves the claim.

Therefore for each 
$j\in \{j_\psi+1,\ldots, h-1\}$ there exists $i_j\in \{1,2\}$ so that 
$\Phi_{i_j}(HN_j(E)\otimes M_{i_j}) \nsubseteq HN_{j}(E)$. 
Then the same argument as in the proof 
of \cite[Prop 3.2]{semistpairs} and  
\cite[Thm. 3.1]{projective} shows that 
\be\label{eq:boundC}
\mu(HN_j(E)/HN_{j-1}(E)) \leq \mu(HN_{j+1}(E)/HN_{j}(E)) 
-\mathrm{deg}(M_{i_j}). 
\ee
Summing inequalities \eqref{eq:boundC} from $j=j_\psi+1$ to $j=h-1$ we obtain 
\[ 
\mu(HN_{j_\psi+1}(E)/HN_{j_\psi}(E)) 
\leq \mu(E/HN_{h-1}(E)) -\sum_{j=j_\psi+1}^{h-1} \mathrm{deg}(M_{i_j}).
\]
Then using inequality \eqref{eq:boundA} we obtain 
\[
\mu_{\mmin}(E_\infty) +\sum_{j=j_\psi+1}^{h-1} \mathrm{deg}(M_{i_j}) 
\leq \mu(E/HN_{h-1}(E)).
\]
which further yields 
\[ 
\mu_{\mmin}(E_\infty) - (h-1)\mmax\{|\mathrm{deg}(M_1)|,|\mathrm{deg}(M_2)|\}
\leq \mu(E/HN_{h-1}(E)).
\]
Since we have established above that $h\leq r$, we finally obtain 
\be\label{eq:boundD} 
- \mu(E/HN_{h-1}(E)) \leq 
-\mu_{\mmin}(E_\infty) +
(r-1)\mmax\{|\mathrm{deg}(M_1)|,|\mathrm{deg}(M_2)|\}
\ee
Taking into account \eqref{eq:boundB}, \eqref{eq:boundD}, inequality  
\eqref{eq:HNineq} implies the existence of the required 
upper bound for $\mu(HN_1(E))=\mu_\mmax(E)$. 

Next suppose $\CE$ is 0-semistable and $\psi$ is identically zero. 
If $\phi$ is nontrivial, boundedness follows from the above argument using 
Lemma (\ref{duallemma}). If $\phi$ is also trivial, definition (\ref{zerostab}) implies that 
the data $(E, \Phi_1, \Phi_2)$ is a semistable Higgs sheaf on 
$X$ as defined in (\ref{Higgs-sheaves}).  Then boundedness follows by a very similar 
argument.

\hfill $\Box$

Lemma (\ref{boundedness}) implies the following corollary by a standard argument. 
\begin{coro}\label{univbound}
The set of isomorphism classes of 
ADHM sheaves of type $(r,e)$ on $X$ which are $\delta$-semistable 
for at least one value $\delta\in \IR_{>0}$ is bounded. 
\end{coro} 

The proof of lemma (\ref{boundedness}) also implies the following. 
\begin{coro}\label{lowerbound} 
Let $r\in \IZ_{\geq 1}$ be a fixed rank. Then there exists $c\in \IZ$ depending only on $r$ so that for 
any $e\in \IZ$, $e<c$ and any  $\delta\in \IR_{>0}$ there are no 
$\delta$-semistable ADHM sheaves 
of type $(r,e)$ on $X$. 
\end{coro} 

{\it Proof.} Suppose $\CE=(E,\Phi_{1,2},\phi,\psi)$ is a $\delta$-semistable ADHM sheaf of type $(r,e)$. If $E$ is semistable, it follows that $\mu(E)\geq \mu_{\mathrm{min}}(E_\infty)$ since there must exist a nontrivial morphism $\psi:E_\infty \to E$. 
If $E$ is not semistable, equation \eqref{eq:boundD} implies that
\[
\mu(E)> \mu_{\mathrm{min}}(E)=\mu(E/HN_{h-1}(E)) \geq
\mu_{\mmin}(E_\infty) -
(r-1)\mmax\{|\mathrm{deg}(M_1)|,|\mathrm{deg}(M_2)|\}.
\]
This proves the claim. 

\hfill $\Box$

 \section{Categorical formulation}\label{catform}

In this subsection we reformulate $\delta$-stability of ADHM sheaves 
as a stability condition in a certain abelian category. 
This will enable us to study the behavior of the moduli spaces 
$\mfm^{ss}_\delta(\CX,r,e)$ as a function of the stability parameter 
$\delta$.  Similar constructions have been carried out for example 
in \cite{Systcoh,King-BN} for moduli spaces of coherent systems. 
The abelian category in question  will be constructed as 
a subcategory of an abelian  category of twisted quiver sheaves on $X$. 
Then the slope stability condition defined below belongs to the class 
of stability conditions studied in \cite{dimred,HKquivers}. 
 
\subsection{An abelian subcategory of ADHM quiver sheaves}\label{absubcat}
Let $X$ be a scheme over a $K$ over $\IC$
 and let $(M_1, M_2)$ be fixed invertible 
sheaves on $X$. Set $M=M_1\otimes_X M_2$ and 
suppose there is a fixed isomorphism $M\simeq K_X^{-1}$ as in the previous 
section. 
ADHM quiver sheaves on $X$ are $(M_1,M_2)$-twisted 
representations of an ADHM quiver in the abelian category 
of quasi-coherent sheaves on $X$. 
Such objects have been considered in the literature in 
\cite{dimred,HKquivers,homquiv,szendroi-2005,modquivers}. 
Basically ADHM quiver sheaves are defined by the same 
data as ADHM sheaves except that the framing data $E_\infty$ 
is not fixed.  More precisely we have 
\begin{defi}\label{quivsheaves}
$(i)$ An ADHM quiver sheaf on $X$ is 
a collection 
$\CE=(E,E_\infty,\Phi_{1,2}, \phi,\psi)$
where $E,E_\infty$ are coherent $\CO_X$-modules 
and 
\[
\Phi_i:E\otimes_X M_i\to E,\qquad \phi:E\otimes_X M\to E_\infty,
\qquad \psi:E_\infty \to E
\]
are morphisms of $\CO_X$-modules satisfying the ADHM relation. 

$(ii)$ A morphism between two ADHM quiver sheaves $\CE$, $\CE'$ is a
pair $(\xi,\xi_\infty)$ of morphisms of $\CO_X$-modules 
\[
\xi:E\to E',\qquad \xi_\infty : E_\infty\to E_\infty'
\]
satisfying the obvious compatibility conditions with the data 
$(\Phi_{1,2},\phi,\psi)$, $(\Phi_{1,2}', \phi', \psi')$.

$(iii)$ An ADHM quiver sheaf $\CE$ will be called locally free if 
$E,E_\infty$ are  locally free $\CO_X$-modules.
\end{defi}

According to \cite{dimred,HKquivers,homquiv,szendroi-2005,modquivers}, ADHM quiver 
sheaves on $X$ form
an abelian category $\CQ_X$. 
Now define a subcategory 
$\CC_X$ of the abelian category ${\CQ}_X$
as follows 
\begin{itemize}
\item[$\bullet$] The objects of $\CC_X$ are coherent ADHM quiver sheaves 
with 
\[
E_\infty = V\otimes \CO_X \]
where $V$ is a finite dimensional vector space over $K$ (possibly 
trivial). We will denote by $v\in \IZ_{\geq 0}$ the 
dimension of $V$.
\item[$\bullet$] Given two objects $\CE,\CE'$ of $\CC_X$ 
a morphism from $\CE$ to $\CE'$ is a morphism $(\xi,\xi_\infty)$ 
of ADHM quiver sheaves so that 
\[
\xi_\infty = f\otimes 1_{\CO_X}
\]
where $f:V \to V'$ is a $K$-linear map.
\end{itemize}

 \begin{rema}\label{isomrema} 
 $(i)$ Note that there is an obvious one-to-one correspondence 
 between ADHM sheaves on $X$ and objects of $\CC_X$  
 with $v=1$. Similarly, there is a obvious one-to-one correspondence between 
 Higgs sheaves, as defined in (\ref{Higgs-sheaves}), and objects of 
 $\CC_X$ with $v=0$. 
 
 $(ii)$ Given two ADHM  sheaves $\CE,\CE'$ on $X$ 
a morphism $(\xi,\lambda):\CE{\to}\CE'$ 
$\CC_{X}$ is a morphism of ADHM sheaves as defined 
in \cite[Def 2.1]{modADHM} if and only if $\lambda=1$. 
This distinction is important in the construction of moduli spaces. 
Note also that if $(\xi,\lambda):\CE{\to}\CE'$ is an morphism in $\CC_{X}$,
with $\lambda\neq 0$,
then $(\lambda^{-1}\xi, 1):\CE{\to}\CE'$ is a morphism of ADHM sheaves. 
\end{rema}

\begin{lemm}\label{abcat}
The category 
 $\CC_{X}$ is abelian. 
\end{lemm}

{\it Proof.} Reduces to a straightforward verification that $\CC_{X}$ contains 
all kernel, images and cokernels of its morphisms, as well as direct sums. 
This is easily done by standard diagram chasing. We will omit the details.

\hfill $\Box$

\subsection{$\delta$-stability in the abelian category}\label{stabsubcat}
In this subsection let $X$ be a smooth projective 
curve over a field $K$ over $\IC$.
Let $\delta \in \IR$. 
Given an object $\CE$ of $\CC_X$, 
the type of $\CE$ is the triple $(r(E),d(E),v)\in 
\IZ_{\geq 0}\times \IZ\times \IZ_{\geq 0}$.
If $r(E)>0$, the $\delta$-slope of $\CE$ is defined by 
\[
\mu_\delta(\CE) = \mu(E)+{v \delta \over r(E)} 
\]
Note that if  $v=0$, $\CE$ is a Higgs sheaf and 
$\mu_{\delta}(\CE)=\mu(\CE)$ is the usual slope of $\CE$ 
for any value of $\delta\in \IR$. Some basic properties of Higgs sheaves and 
their moduli are summarized for convenience in Appendix (\ref{higgsapp}). 

Let $\delta \in \IR$. 
Then we define $\delta$-(semi)stability for objects of 
$\CC_{X}$ as follows. 
\begin{defi}\label{abcatstab}
An object $\CE$ of $\CC_{X}$ of type $(r,e,v)\in \IZ_{\geq 0}\times 
\IZ \times \IZ_{\geq 0}$ will be called $\delta$-(semi)stable if the 
following inequality holds  for any proper nontrivial subobject $0\subset \CE'\subset \CE$ in 
$\CC_X$ 
of type $(r',e',v')$ 
\be\label{eq:abstabcond}
 r\left(e'+v'\delta\right)\ (\leq)\
 r'\left(e+ v \delta\right)
 \ee
\end{defi}
 
 \begin{rema}\label{rankzeroobj}
 $(i)$ Note that the object $O = (0, \IC, 0,0,0,0)$ is stable according to 
 definition (\ref{abcatstab}) for any value of $\delta \in \IR$. 
 
 $(ii)$ Note that if $v(\CE)=0$, the $\delta$-stability condition reduces to the 
 standard slope stability condition for Higgs sheaves. 
 \end{rema}
 
 The following results follow by standard manipulations of the stability 
 conditions. Details will be omitted. 
  \begin{lemm}\label{endlemma}
 Let $\CE$ be a 
 $\delta$-semistable object of $\CC_{X}$ of rank $r\in \IZ_{\geq 1}$ 
 for some $\delta \in \IR_{>0}$. 
 Then the following hold. 
 
 $(i)$ $E$ is a locally free $\CO_X$-module. Moreover, if in addition 
 $v(\CE)>0$, the morphism $\psi:\CO_X\to E$ must be nontrivial. 
 
 $(ii)$ Any nontrivial endomorphism $(\xi,f):\CE\to \CE$ 
 in $\CC_X$ must be an isomorphism. If in addition 
 the ground field $K$ is algebraically closed, 
 the endomorphism ring of $\CE$ is isomorphic to $K$. 
 \end{lemm}
   
\begin{lemm}\label{equivstab}
Let $\delta\in \IR$. Let 
$\CE$ be a locally free object of $\CC_{X}$ of rank $r(\CE)\geq 1$ and $v(\CE)=1$.
Then $\CE$ is $\delta$-(semi)stable if and only if it is 
$\delta$-(semi)stable as an ADHM sheaf.
\end{lemm}

Given remark (\ref{isomrema}), an immediate consequence of lemmas 
 (\ref{endlemma}), (\ref{equivstab}) is 
 \begin{coro}\label{automADHM} 
 Suppose $K$ is algebraically closed and let $\CE$ be a nontrivial $\delta$-stable 
 ADHM sheaf on $X$ for some 
 $\delta\in \IR_{>0}$. Then the automorphism group  of $\CE$ is trivial. 
 \end{coro}

To conclude this subsection, note that since the category $\CC_{X}$ is 
noetherian and 
artinian the standard properties of $\delta$-(semi)stable 
objects hold. That is we have 
\begin{prop}\label{abstabprop}
$(i)$ Harder-Narasimhan filtrations in $\CC_{X}$ exist and satisfy 
the same properties 
as Harder-Narasimhan filtrations of coherent sheaves on smooth 
projective varieties. 

$(ii)$ The subcategory of $\delta$-(semi)stable objets of $\CC_{X}$ 
with fixed $\delta$-slope is noetherian and artinian. The simple objects in this 
subcategory are precisely the $\delta$-stable objects. 
In particular Jordan-H\"older filtration exist and satisfy the same properties 
as Jordan-H\"older filtrations of semistable coherent sheaves on smooth 
projective varieties. 
\end{prop}

\subsection{Extensions of ADHM quiver sheaves}\label{extsect}
Next we prove some basic homological algebra results
for ADHM quiver sheaves.  The homological algebra of 
quiver sheaves without relations has been treated in detail in 
\cite{homquiv}. 
Our task is to generalize some of the results of \cite{homquiv}
to ADHM quiver sheaves, which have quadratic relations. 
In this subsection we will take $X$ to be a separated scheme 
of finite type over a  field $K$  over $\IC$
and will employ $\check{\mathrm{C}}$ech 
resolutions rather than injective resolutions as in \cite{homquiv}. 
Let $M_1,M_2$ be fixed invertible sheaves on $X$. 

The main result of this section is 
\begin{prop}\label{quivextA}
Let $\CE'=(E',E'_\infty,\Phi'_{1,2},\phi',\psi')$, 
$\CE''=(E'',E''_\infty,\Phi''_{1,2},\phi'',\psi'')$ be coherent 
locally free ADHM quiver sheaves on $X$. 
Consider the following complex $\CC(\CE'',\CE')$ of coherent 
locally free $\CO_X$-modules 
\be\label{eq:hypercohA} 
\begin{aligned} 
0 \to \begin{array}{c} \lochom_{X}(E'',E') \\ \oplus \\ 
\lochom_{X}(E''_{\infty},E'_{\infty}) \\ \end{array}
& {\buildrel d_1\over \longto} 
\begin{array}{c}  \lochom_{X}(E''\otimes_{X}M_1 ,E') \\ \oplus \\
\lochom_{X}(E''\otimes_{X} M_2, E') \\ \oplus \\ 
\lochom_{X}(E''\otimes_{X} M,E_{\infty}')\\ \oplus \\ 
\lochom_{X}(E''_{\infty},E') \\ \end{array} 
 {\buildrel d_2\over \longto} 
\lochom_{X}(E''\otimes_{X}M,E') \to 0 \\
\end{aligned}
\ee
where 
\[
\begin{aligned}
d_1(\alpha,\alpha_\infty) = 
(& -\alpha \circ \Phi_{1}'' +\Phi_{1}'\circ (\alpha\otimes 1_{M_1}), 
-\alpha \circ \Phi_{2}''+\Phi_{2}'\circ (\alpha\otimes 1_{M_2}),\\ 
& -\alpha_\infty\circ \phi'' +\phi' \circ (\alpha
\otimes 1_M), 
-\alpha\circ \psi'' +\psi' \circ \alpha_\infty)\\
\end{aligned}  
\]
for any local sections $(\alpha,\alpha_\infty)$ of 
the first term 
and 
\[
\begin{aligned}
d_2(\beta_1,\beta_2,\gamma, \delta) = &
\beta_1 \circ (\Phi''_2\otimes 1_{M_1}) -
\Phi_{2}'\circ (\beta_1\otimes 1_{M_2}) 
- \beta_2\circ (\Phi''_{1}\otimes 1_{M_2})\\
 & + \Phi_{1}'\circ (\beta_2\otimes 1_{M_1}) + 
\psi'\circ \gamma + \delta \circ \phi''\\
\end{aligned}
\]
for any local sections $(\beta_1,\beta_2,\gamma, \delta)$ 
of the middle term. The degrees of the three terms in 
\eqref{eq:hypercohA} are $0,1,2$ respectively. 

Then there are group isomorphisms 
\be\label{eq:adhmquivext} 
\mathrm{Ext}^k_{{\mathcal Q}_X}(\CE'',\CE')\simeq \IH^k(X, \CC(\CE'',\CE'))
\ee
for $k=0,1$, where $\mathrm{Ext}^k_{{\mathcal Q}_X}(\CE'',\CE')$ 
denote extension groups in the abelian category of ADHM quiver sheaves on $X$. 
\end{prop}

{\it Proof.}
Since $X$ is a separated scheme of finite type over $\IC$, it 
admits affine open covers and we can employ $\check{\mathrm{C}}$ech
resolutions in the construction of the hypercohomology double complex
associated to $\CC(\CE'',\CE')$.  
Then the correspondence stated in lemma (\ref{quivextA}) 
follows by repeating the proof of \cite[Prop. 4.5]{modADHM} 
based on $\check{\mathrm{C}}$ech  
cochain computations in the present context. 
We will omit the details. 

\hfill $\Box$

Proposition (\ref{quivextA}) implies 
\begin{coro}\label{Cextcoro}
Suppose $X$ is connected and proper over $K$. Let $\CE', \CE''$ be locally free 
objects of $\CC_X$ 
with $v(\CE')+v(\CE'') \leq 1$. Then there are group isomorphisms 
\be\label{eq:CextA}
\mathrm{Ext}^k_{{\mathcal C}_X}(\CE'',\CE')\simeq \IH^k(X, \CC(\CE'',\CE'))
\ee
for $k=0,1$, where $\mathrm{Ext}^k_{{\mathcal C}_X}(\CE'',\CE')$ 
denote extension groups in the abelian category $\CC_X$. 
\end{coro}

{\it Proof.} 
The case $k=0$ follows from the fact that $\CC_X$ is a full subcategory of 
${\mathcal Q}_X$ if $X$ is connected and proper over $K$.
For $k=1$ note that there is a natural injective homomorphism 
\[
\mathrm{Ext}^1_{{\mathcal C}_X}(\CE'',\CE')\hookrightarrow 
\mathrm{Ext}^1_{{\mathcal Q}_X}(\CE'',\CE').
\]
If $v(\CE')+v(\CE'')\leq 1$, it follows that at least one of $v(\CE')$, $v(\CE'')$ 
vanishes. Then it is straightforward to check that the above homomorphism 
is also surjective. 

\hfill $\Box$

\subsection{Stability conditions and moduli spaces of ADHM quiver sheaves}\label{adhmquivsect} 
Several results concerning stability conditions and moduli spaces 
of ADHM quiver sheaves are summarized in this section
for future reference. Here $X$ will be a smooth projective curve over $\IC$. 
Note that very general stability conditions for quiver sheaves 
have been introduced in \cite{HKquivers}. Specializing these conditions to objects of ${\CQ_X}$ 
results in the following. 

\begin{defi}\label{adhmquivstab} 
Let $(\sigma, \tau)\in \IQ_{>0}\times \IQ$ be fixed stability parameters. 
An object $\CE$ of ${\CQ}_X$ is $(\sigma,\tau)$-semistable 
if any proper nontrivial quiver subsheaf $0\subset \CE'\subset \CE$ 
satisfies 
\be\label{eq:quivsheafstab}
(r(E) + \sigma r(E_\infty))(d(E')+ \tau r(E_\infty')) \ (\leq) \  
  (r(E') + \sigma r(E'_\infty))(d(E)+ \tau r(E_\infty))
  \ee
  \end{defi}
 Note that $(\sigma,\tau)$-semistable ADHM quiver sheaves must be locally free, 
 and it is sufficient to test the stability condition on saturated subobjects. 
 Moreover, the automorphism group of any stable ADHM quiver sheaf 
 is canonically isomorphic to $\IC^\times$. 
 
 Moduli problems for general quiver sheaves have been treated in \cite{modquivers}. In particular, \cite[Thm 3.6.1]{modquivers} 
 proves that there exists a coarse quasi-projective moduli 
 space of quiver sheaves on a smooth projective variety 
 subject to a stability condition involving weighted flags of subobjects. 
 Then \cite[Thm. 3.7.1]{modquivers} proves that this stability 
 condition stabilizes to a conventional one of the type studied 
 in \cite{HKquivers} if a certain stability parameter is taken very large. 
 Moreover, \cite[Thm 3.7.2]{modquivers} establishes the existence 
 of a proper Hitchin morphism for moduli spaces of quiver sheaves 
 defined by polynomial invariants. In conclusion, in the present 
 case, \cite[Thm. 3.6.1, Thm. 3.7.1, Thm 3.7.2]{modquivers} imply the following result. 
 \begin{theo}\label{adhmquivmod} 
 For fixed numerical invariants $r,r_\infty\in \IZ_{\geq 1}$, $e,e_\infty\in \IZ$ there exists a
quasi-projective 
coarse moduli scheme ${\mathcal Q}^{ss}_{(\sigma,\tau)}(\CX,r,e; r_\infty,e_\infty)$ 
parameterizing S-equivalence classes 
 of $(\sigma,\tau)$-semistable ADHM quiver sheaves with fixed numerical invariants 
$r,r_\infty\in \IZ_{\geq 1}$, $e,e_\infty\in \IZ$. This moduli scheme contains an open subscheme 
${\mathcal Q}^{s}_{(\sigma,\tau)}(\CX,r,e; r_\infty,e_\infty)$ parameterizing isomorphism 
classes of $(\sigma,\tau)$-stable objects. 

Moreover, there exists a proper generalized Hitchin morphism 
$h_{Q}: {\mathcal Q}^{ss}_{(\sigma,\tau)}(\CX,r,e; r_\infty,e_\infty)\to \IV$, 
with $\IV$ an affine space defined by the polynomial invariants of the morphism 
data $(\Phi_{1,2},\phi,\psi)$ of an ADHM quiver sheaf $\CE=(E,E_\infty, \Phi_{1,2}, 
\phi,\psi)$. 
\end{theo}

\begin{rema}\label{GITconstr} 
Note that the moduli space ${\mathcal Q}^{ss}_{(\sigma,\tau)}(\CX,r,e; r_\infty,e_\infty)$ 
is constructed in \cite{modquivers}  as a GIT 
quotient of a parameter space $P$ by an affine reductive group $G$. In particular 
$P$ is equipped with a G-linearized invertible sheaf $L$. The data 
$(P,G,L)$ determines an ample line bundle $\CL$ on the moduli space.
This construction also implies that the quotient stack $[P/G]$ is isomorphic to the 
groupoid ${\mathfrak Q}^{ss}_{(\sigma,\tau)}(\CX,r,e; r_\infty,e_\infty)$ 
of flat families of  $(\sigma,\tau)$-semistable ADHM quiver sheaves on $X$. 
Moreover by analogy with \cite[Ex. 7.7]{goodmoduli} 
there is a natural morphism
\be\label{eq:goodmorphism} 
{\mathfrak q}: {\mathfrak Q}^{ss}_{(\sigma,\tau)}(\CX,r,e; r_\infty,e_\infty)\to 
{\mathcal Q}^{ss}_{(\sigma,\tau)}(\CX,r,e; r_\infty,e_\infty)
\ee
which satisfies the conditions of \cite[Thm. 4.14]{goodmoduli}. In particular 
${\mathfrak q}$ is universally closed.  Moreover the restriction 
of ${\mathfrak q}$ to the substack of $(\sigma,\tau)$-stable objects is a
$\IC^\times$-gerbe over ${\mathcal Q}^{s}_{(\sigma,\tau)}(\CX,r,e; r_\infty,e_\infty)$.
\end{rema}
 
Finally, 
note that there is a relation between $(\sigma,\tau)$-stability in the abelian category $\CQ_X$ and 
$\delta$-stability in the subcategory 
$\CC_X$. This will be explained in detail below because it will be used  in the proof of Theorem (\ref{modthm}).

\begin{lemm}\label{sigmataudelta}
Fix $(r,e)\in \IZ_{\geq 1}\times \IZ$, $(r_\infty,e_\infty)=(1,0)$ and 
$\delta\in \IQ_{>0}$. Then 
there exists $\sigma_0\in \IQ_{>0}$ (depending on $(r,e,\delta)$), such that 
the following hold
for any $0<\sigma < \sigma_0$. 

$(i)$ If $\CE$ is an 
object of $\CC_X$ which is $(\sigma, \delta)$-stable as an object of 
$\CQ_X$, then it is $\delta$-semistable. 

$(ii)$ If $\CE$ is a $\delta$-stable object of $\CC_X$ of type $(r,e)$, then 
$\CE$ is $(\sigma,\delta)$-stable as an object of $\CQ_X$.
 \end{lemm}
 
 {\it Proof.} 
 First note that the set of isomorphism classes of objects of $\CC_X$ 
  with fixed numerical invariants as above which 
 are $(\sigma,\delta)$-semistable for some $\sigma \in \IQ_{>0}$ is bounded. 
 The proof is very similar to the proof of Lemma (\ref{boundedness}). 
 
 In order to prove $(i)$ suppose for any $\sigma \in \IQ_{>0}$
 there exists a $(\sigma, \delta)$-stable object $\CE$ of $\CQ_X$ of type $(r,e,1,0)$ which 
 is not $\delta$-semistable. 
 Given such an object $\CE$ there exists a proper nontrivial subobject $0\subset \CE'\subset \CE$ 
 in $\CC_X$ such that 
 \[
 r(e'+v'\delta) > r'(e+\delta). 
 \] 
 Moreover, $E'\subset E$ may be assumed saturated, and $0<r'<r$. 
 Since $\CE'$ is also a subobject of $\CE$ in $\CQ_X$, 
 \[
 (r+\sigma) (r'+e' \delta) < (r'+v'\sigma)(e+\delta). 
 \]
  If $v'=0$, the above inequalities imply 
  \[
  {r'\over r}(e+\delta) < e' < {r'\over r+\sigma}(e+\delta) 
  \] 
 Since $\sigma>0$, this cannot hold unless $e+\delta<0$. 
 
 If $v'=1$, it follows that 
 \[
 {r'\over r}(e+\delta) < e' +\delta< {r'+\sigma \over r+\sigma}(e+\delta) 
 \]
 which cannot hold unless $e+\delta>0$. 
 
 In both cases, it follows that the absolute values of  $|e'|$ 
is bounded by a constant depending only on $(r,e, \delta)$ 
  since $0< r'<r$. Then \cite[Lemm. 1.7.9]{huylehn} implies 
   that the set of isomorphism classes 
  of all such destabilizing subobjects $\CE'$ is bounded since the set of all 
  $(\sigma,\tau)$-semistable objects $\CE$  for some value of $\sigma>0$ is bounded. 
Therefore the slope $e'/r'$ can only take values in a finite subset of $\IQ$, which is 
independent of $\sigma$. This leads to a contradiction for sufficiently small $\sigma$ 

The proof of $(ii)$ is similar. 

\hfill $\Box$
  
\section{Chamber structure}\label{chambersect}

The goal of this section is to study the behavior of the 
$\delta$-stability condition on the parameter $\delta\in \IR_{>0}$ 
keeping the data $\CX=(X,M_1,M_2)$, $E_\infty =\CO_X$ as 
well as the type $(r,e)\in \IZ_{\geq 1}\times \IZ$ 
fixed. The ground field will be a field $K$ over $\IC$ as in section (\ref{adhmstab}). 

\subsection{Critical stability parameters}\label{crtstabpar}
This section establishes the existence of a chamber 
structure of the positive real axis so that the set of $\delta$-stable 
ADHM sheaves is constant in each chamber.

First consider the case of ADHM sheaves of rank $r=1$. 
\begin{lemm}\label{rankone}
Let $\CE=(E,\Phi_1,\Phi_2,\phi,\psi)$ be a locally free 
ADHM sheaf of type $(1,e)$, $e\in \IZ$, on 
$X$ so that $\psi$ is nontrivial. Then $\CE$ is 
$\delta$-stable for any stability parameter $\delta\in \IR_{>0}$. 
In particular the moduli space $\mfm_{\delta}(\CX,1,e)$ is independent of 
$\delta \in \IR_{>0}$.
\end{lemm}

{\it Proof.} 
Since $r=1$, $E$ has no nontrivial proper saturated subsheaves 
$0\subset E'\subset E$. Hence, since $\psi$ is nontrivial 
the $\delta$-stability conditions are automatically satisfied 
for any $\delta\in \IR_{>0}$. 

\hfill $\Box$

\begin{rema}\label{rankonerema} 
Given lemma (\ref{rankone}), a rank one locally free ADHM sheaf on $X$ 
with $\psi\neq 0$ will be called in the following stable, without any reference 
to a stability parameter.  
\end{rema}

Next let $r\geq 2$. 
Let $\delta \in \IR_{>0}$ be stability parameter. Suppose there exists 
a $\delta$-semistable ADHM sheaf $\CE$ of type $(r,e)$ on $X$ 
which is not $\delta$-stable. 
Then definition (\ref{deltastable}) implies that $\delta$ must be of the form 
\be\label{eq:crtform}
\delta = {re'-er'\over r'} \qquad \mathrm{or} \qquad 
\delta = {er'-re'\over r-r'} 
\ee
for some $1\leq r'\leq r-1$, $e'\in \IZ$. 

\begin{defi}\label{crtval} 
A stability parameter $\delta\in \IR_{>0}$ is called numerically critical 
of type $(r,e)\in \IZ_{\geq 2}\times \IZ$ if it is of the form \eqref{eq:crtform}. 
\end{defi} 

Let $\Delta_{(r,e)}\subset \IR_{>0}$ denote the set of numerically critical parameter of fixed 
type $(r,e)\in \IZ_{\geq 2}\times \IZ$. 
Since all such parameters are rational numbers with denominators 
in the finite set $\{1,,\ldots, r-1\}$ it follows that there exists 
an isomorphism $\Delta_{(r,e)}\simeq \IZ_{>0}$. For each $n\in 
\IZ_{>0}$ let $\delta_n \in \Delta_{(r,e)}$ denote the corresponding numerically
critical parameter. We will also set $\delta_0=0$ in order to simplify the 
exposition. Then the following result holds by standard arguments. 
\begin{lemm}\label{chamber}
Let $(r,e)\in \IZ_{\geq 2}\times \IZ$ be a fixed type. Then the following hold 
\begin{itemize}
\item[$(i)$]  
For any $n\in \IZ_{\geq 0}$ and any 
$\delta\in (\delta_{n}, \delta_{n+1})$
an ADHM sheaf of type $(r,e)$ is $\delta$-semistable 
if and only if it is $\delta$-stable.
\item[$(ii)$] For any $n\in \IZ_{\geq 0}$, the 
set of $\delta$-stable ADHM sheaves
is constant for $\delta \in (\delta_{n}, \delta_{n+1})$ 
\end{itemize} 
\end{lemm}

\begin{defi}\label{asympstab}
A locally free ADHM sheaf $\CE$ on $X$ is asymptotically stable if 
$\psi$ is not identically zero and there is no nontrivial $\Phi$-invariant
proper saturated subsheaf 
$0\subset E'\subset E$ so that $\mathrm{Im}(\psi) \subseteq E'$.
\end{defi}

Then the following boundedness result holds. Since proof is analogous to the 
proof of Lemma (\ref{boundedness}) the  details will be omitted. 

\begin{lemm}\label{asympbound} 
The set of isomorphism classes of locally free coherent sheaves $E$ 
of type $(r,e)\in \IZ_{\geq 2}\times \IZ$ 
on $X$ with the property that $E$ is the underlying 
sheaf of an asymptotically stable ADHM sheaf $\CE$ is bounded.
\end{lemm}

\begin{lemm}\label{largedelta}
Let $(r,e)\in \IZ_{\geq 2}\times \IZ$ be fixed as above. 
Then there exists  $\delta_\infty\in \IR_{>0}$ depending only on 
$(r,e)$ 
so that the following statements hold  
for any $\delta > \delta_\infty$
\begin{itemize}
\item[$(i)$] 
An ADHM sheaf $\CE$ of type $(r,e)$ on $X$ is 
$\delta$-semistable if and only if it is asymptotically stable. 
\item[$(ii)$]
An ADHM sheaf $\CE$ of type $(r,e)$ on $X$ is 
$\delta$-semistable if and only if it is $\delta$-stable
\end{itemize}
\end{lemm} 

{\it Proof.} 
First we prove that there exists $c_1\in \IR_{>0}$ depending only 
on $(r,e)$ so that  any 
$\delta$-semistable ADHM sheaf on $X$ of type $(r,e)$ with $\delta>c_1$ 
is asymptotically stable. 
Recall \cite[Lemma 2.4]{modADHM} 
that given any ADHM sheaf $\CE=(E,\Phi_{1,2},\phi,\psi)$ 
with $\psi\neq 0$ 
there is a canonical nontrivial $\Phi$-invariant saturated subsheaf $E_0\subset E$ 
so that 
$\mathrm{Im}(\psi)\subset E_0$. 
Moreover, by construction $E_0$ is a subsheaf of 
any $\Phi$-invariant saturated subsheaf $E'\subset E$ so that 
$\mathrm{Im}(\psi)\subset E'$. 

Since $E_0$ is canonically constructed 
in terms of the data $\CE$, it follows that the set of the
 sheaves $E_0$ associated to all $\delta$-semistable 
 ADHM sheaves $\CE$ of type $(r,e)$, with $\delta\in \IR_{>0}$, is bounded. 
 Therefore the numerical invariants $(r(E_0),d(E_0))$ can take only a finite set of values. 
  This implies that there exists $c_1\in \IR_{>0}$ so that for all 
 $\delta> c_1$ we have 
 \be\label{eq:destabcond}
 \mu(E_0)+{\delta\over r(E_0)}  > \mu_{\delta}(\CE) 
\ee
whenever $E_0$ is a proper subsheaf of $E$. 

Suppose there exists a
$\delta$-semistable ADHM sheaf $\CE$, $\delta >c_1$, which is not asymptotically 
stable. Therefore there exists a nontrivial $\Phi$-invariant proper saturated 
subsheaf $0\subset E'\subset E$ so that $\mathrm{Im}(\psi)\subseteq E'$. 
As observed above, 
by construction $E_0$ must be a subsheaf of 
any such subsheaf, hence in particular $E_0$ 
will be  a proper subsheaf of $E$. 
This yields a contradiction since then 
inequality \eqref{eq:destabcond} implies that $\CE$ is not 
$\delta$-semistable. 

Next we prove that there exists $c_2\in \IR_{>0}$ so that 
 any asymptotically stable ADHM sheaf 
$\CE$  is $\delta$-stable for all $\delta>c_2$. 
According to lemma (\ref{asympbound}), the set of 
isomorphism classes of locally free sheaves $E$ of type $(r,e)$ so that 
$E$ is the underlying sheaf of an asymptotically stable ADHM sheaf 
is bounded. 
This implies that there exists a positive 
constant $C'\in \IR_{>0}$ depending only on 
$(r,e)$ so that $\mu_{\mathrm{max}}(E) < C'$ 
for any such locally free sheaf $E$. 
It follows that there exists $c_2\in \IR_{>0}$ so that 
for any $\delta>c_2$, condition \eqref{eq:slopecondB}
 of definition (\ref{deltastable}) 
is automatically satisfied for any asymptotically stable ADHM sheaf $\CE$, 
and any $\Phi$-invariant nontrivial proper saturated subsheaf 
$0\subset E'\subset E$ so that $E'\otimes_X M \subset \mathrm{Ker}(\phi)$.
Since the stability condition $(ii)$ of definition (\ref{deltastable}) is trivially satisfied for 
asymptotically stable ADHM sheaves, the claim follows. 

In order to conclude the proof of proposition (\ref{largedelta}), take $\delta_\infty = 
\mathrm{max}\{c_1,c_2\}$. 

\hfill $\Box$

\begin{lemm}\label{finitecrt}
Let $(r,e)\in \IZ_{\geq 2}\times \IZ$ be a fixed type. Then 
there exists $\delta_N \in \Delta_{(r,e)}\cup \{0\}$, $N\geq 0$, 
depending 
only on $(r,e)$ so that for all $\delta>\delta_N$, 
any $\delta$-semistable ADHM sheaf 
of type $(r,e)$ is asymptotically stable. 
\end{lemm} 

{\it Proof.} 
Follows directly from lemmas (\ref{chamber}) and (\ref{largedelta}). 

\hfill $\Box$

\begin{defi}\label{crtvaldef}
Suppose the integer $N$ found in lemma 
(\ref{finitecrt}) is nonzero. Then the stability 
parameters $\delta_i\in \IR_{>0}$, $i=1,\ldots, N$ 
 will be called critical values of type $(r,e)$. 
 Moreover, a parameter $\delta \in \IR_{>0}$ will be called noncritical of type 
 $(r,e)$ if $\delta \notin \{\delta_1,\ldots, \delta_N\}$. 
\end{defi} 

\subsection{Wallcrossing behavior}\label{wallsect}
This subsection analyzes  the behavior of $\delta$-stable ADHM sheaves 
as $\delta$ specializes to a critical value. In order to simply the exposition we 
 will formally set $\delta_0 =0$ as above, and $\delta_{N+1}=+\infty$. 
\begin{lemm}\label{specstab} 
Let $(r,e)\in \IZ_{\geq 2}\times \IZ$ be a fixed type so that the integer $N$ 
in proposition (\ref{chamber}) is nonzero. Then the following hold 

$(i)$ Let $\CE$ be a $\delta$-stable ADHM sheaf on $X$ of type $(r,e)$, 
with $\delta \in 
(\delta_i,\ \delta_{i+1})$ for some $i=1,\ldots, N$.
Then $\CE$ is $\delta_i$-semistable 
and it has a Jordan-H\"older filtration in the abelian category $\CC_0$ 
of the form 
\be\label{eq:JHfiltrA}
0=JH_0(\CE) \subset JH_1(\CE) \subset JH_2(\CE) \subset\cdots \subset JH_{j-1}(\CE) 
\subset JH_j(\CE)=\CE, \quad j\geq 1
\ee
so that $r(JH_l(\CE))\geq 1$ for any $1\leq l\leq j$ and $r(JH_l(\CE))<r$,
$v(JH_l(\CE)) = 0$ for $0\leq l\leq j-1$. 
In particular $\CE$ is either $\delta_i$-stable or 
there is a nontrivial extension  in the abelian 
category $\CC_X$ of the form 
\be\label{eq:extensionA}
0\to \CE' \to \CE\to \CE''\to 0
\ee
where $\CE'$ is a semistable Higgs sheaf of rank $r(\CE')\geq 1$
and $\CE''$ is a $\delta_i$-stable 
ADHM sheaf of rank $r(\CE'')\geq 1$ and $\mu_{\delta_i}(\CE') = \mu_{\delta_i}(\CE'')
=\mu_{\delta_i}(\CE)$. 

$(ii)$ Let $\CE$ be a $\delta$-stable ADHM sheaf on $X$ of type $(r,e)$, 
with $\delta \in 
(\delta_{i-1},\ \delta_{i})$ for some $i=1,\ldots, N$.
Then $\CE$ is $\delta_i$-semistable 
and it has a Jordan-H\"older filtration in the abelian category $\CC_X$ 
of the form \eqref{eq:JHfiltrA} 
where $r(JH_l(\CE))\geq 1$, $v(JH_l(\CE)) = 1$ for all $1\leq l\leq j$, 
and 
$r(JH_l(\CE))<r$ for all $0\leq l\leq j-1$. 
In particular $\CE$ is $\delta_i$-stable
or there is a  nontrivial extension in the abelian 
category $\CC_X$ of the form 
\be\label{eq:extensionB}
0\to \CE' \to \CE\to \CE''\to 0
\ee
where $\CE'$ is a $\delta_i$-stable ADHM sheaf of rank $r(\CE')\geq 1$, 
$\CE''$ is a semistable Higgs sheaf of rank $r(\CE'') \geq 1$, and 
$\mu_{\delta_i}(\CE') = \mu_{\delta_i}(\CE'')=\mu_{\delta_i}(\CE)$.

$(iii)$ Let $\CE$ be a $\delta_i$-stable ADHM sheaf on $X$ of type $(r,e)$ for some 
$i=1,\ldots,N$. Then $\CE$ is 
$\delta$-stable for any $\delta \in (\delta_{i-1},\, \delta_{i+1})$.
\end{lemm}

{\it Proof.} 
Since the proofs are very similar, details will be provided only for the first statement. 

Lemma (\ref{chamber}) implies that 
$\CE$ is $\gamma$-stable 
for any $\gamma \in (\delta_i,\ \delta_{i+1})$. Then any $\Phi$-invariant 
nontrivial proper saturated subsheaf $\mathrm{Im}(\psi)\subseteq
 E'\subset E$ must satisfy 
\be\label{eq:gammacondA}
\mu(E') +{\gamma \over r(E')}
< \mu(E) + {\gamma \over r}
\ee
for any $\gamma \in (\delta_i, \delta_{i+1})$. 
Similarly, any $\Phi$-invariant 
nontrivial proper saturated subsheaf $0 \subset
 E'\otimes_X M \subseteq\mathrm{Ker}(\phi)$ must satisfy 
\be\label{eq:gammacondB}
\mu(E') < \mu(E) + {\gamma \over r}
\ee
for any $\gamma \in (\delta_i, \delta_{i+1})$. 

In the first case, it follows that 
\[
\mu(E') +{\delta_i \over r(E')}
\leq  \mu(E) + {\delta_i \over r}
\]
since $r(E') < r$. 

In the second case, suppose there exists such a subsheaf $E'$ so that 
\[
\mu(E') > \mu(E)+{\delta_i\over r}
\]
Then Grothendieck's lemma \cite[Lemma 1.7.9]{huylehn} 
implies that for fixed $E$ the set of isomorphism classes of such subsheaves is bounded. 
Then it follows that there exists $\gamma \in (\delta_i, \delta_{i+1})$ 
so that 
\[
\mu(E') > \mu(E) +{\gamma\over r}.
\]
This would contradict $\gamma$-stability. Therefore $\CE$ must be 
$\delta_i$-semistable, and it has a Jordan-H\"older filtration 
of the form \eqref{eq:JHfiltrA} according to proposition (\ref{abstabprop}).
It is straightforward to check that none of the objects $JH_l(\CE)$, $1\leq l \leq j$ 
may have rank zero and none of the objects $JH_l(\CE)$, $0\leq l \leq j-1$
may have rank equal to $r$. 
If the length of the filtration is $j=1$, it follows that $\CE$ is $\delta_i$-stable. 

Suppose $j\geq 2$. 
By the general properties of Jordan-H\"older filtrations, all subobjects
$JH_l(\CE)\subseteq \CE$, $l=1,\ldots,j$ must have the same 
$\delta_i$-slope as $\CE$. Let us denote by $E_l$ the underlying 
locally free $\CO_X$-module of the ADHM sheaf $JH_l(\CE)$, 
$l=1,\ldots, j$. Let $v_l = v(JH(\CE_l))$, $l=1,\ldots,j$. Then  
\[
\mu(E_l) +{v_l\delta_i \over r(E_l)} 
= \mu(E) + {\delta_i\over r}
\]
for all $l=1,\ldots, j$. 
However, since $\CE$ is $\delta$-stable,
\[
\mu(E_l) + {v_l\delta \over r(E_l)} 
< \mu(E) + {\delta\over r}
\]
must also hold for all $l=1,\ldots, j-1$. These inequalities imply 
\[ 
{v_l(\delta -\delta_i) \over r(E_l)} < 
{\delta -\delta_i\over r}
\]
for all $l=1,\ldots, j-1$. Since $r(E_l)<r$ for $l\neq j$, it follows that 
$v_l=0$ for all $l=1,\ldots, j-1$.  This implies that the last quotient 
$JH_j(\CE)/JH_{j-1}(\CE)$ is a $\delta_i$-stable ADHM sheaf on $X$.
Then the exact sequence \eqref{eq:extensionA} is obtained by 
setting $\CE'=JH_{j-1}(\CE)$, $\CE''=JH_{j}(\CE)$. 
The extension \eqref{eq:extensionA} 
has to be nontrivial because $\CE$ is $\delta$-stable, hence 
indecomposable. 

\hfill $\Box$

The following is an immediate consequence of   lemma (\ref{specstab}).

\begin{coro}\label{wallHN} 
Under the assumptions of lemma (\ref{specstab}), 
let $\delta_i$, $i=1,\ldots, N$ be a critical stability parameter of 
type $(r,e)$ and let $\delta_-\in (\delta_{i-1},\ \delta_i)$, 
$\delta_+\in (\delta_i, \ \delta_{i+1})$ be noncritical 
stability parameters. 

$(i)$ Suppose $\CE$ is a $\delta_+$-stable  ADHM sheaf on $X$ of type $(r,e)$ 
which is not $\delta_-$-stable. Then $\CE$ is strictly $\delta_i$-semistable, 
and in particular it fits in a nontrivial extension of the form \eqref{eq:extensionA}.
Moreover, the one step 
filtration 
$0\subset \CE'\subset \CE$ determined by \eqref{eq:extensionA} 
is a Harder-Narasimhan filtration for $\CE$ with respect with 
$\delta_-$-stability.

$(i)$ Suppose $\CE$ is a $\delta_-$-stable  ADHM sheaf on $X$ of type $(r,e)$ 
which is not $\delta_+$-stable. Then $\CE$ is strictly $\delta_i$-semistable, 
and in particular it fits in a nontrivial extension of the form \eqref{eq:extensionB}.
Moreover, the one step 
filtration 
$0\subset \CE'\subset \CE$ determined by \eqref{eq:extensionB} 
is a Harder-Narasimhan filtration for $\CE$ with respect with 
$\delta_+$-stability.
\end{coro}

For future reference, let us record the following partial converse to lemma 
(\ref{specstab}). 
\begin{lemm}\label{minusdestabobj} 
Under the assumptions of lemma (\ref{specstab}) 
let $\delta_i\in \IR_{>0}$ be a critical stability parameter of 
type $(r,e)\in \IZ_{\geq 2}\times \IZ$. Then the following hold. 

$(i)$ There exists $0<\epsilon_+ <\delta_{i+1}-\delta_i$, 
so that the following holds for any $\delta_+\in (\delta_i, \ \delta_i+\epsilon_+)$. 
A locally free ADHM sheaf $\CE$ of type $(r,e)$ 
on $X$ is $\delta_i$-semistable if and only if 
it is either $\delta_+$-stable or there exists a unique filtration of the form 
\be\label{eq:plusfiltr} 
0\subset \CE'\subset \CE 
\ee
so that $\CE'$ is a $\delta_+$-stable ADHM sheaf of rank $r(\CE')\geq 1$, 
and $\CE''=\CE/\CE'$
is a semistable Higgs sheaf of rank $r(\CE'')\geq 1$  satisfying
\be\label{eq:hnslopecondA}
\mu_{\delta_+}(\CE') > \mu(\CE'') \qquad 
\mu_{\delta_i}(\CE') = \mu(\CE'').
\ee

$(ii)$ There exists $0< \epsilon_-< \delta_i-\delta_{i-1}$ 
so that the following holds for any  
$\delta_-\in (\delta_i-\epsilon_-, \ \delta_i)$. 
A locally free ADHM sheaf $\CE$ of type $(r,e)$ 
on $X$ is $\delta_i$-semistable if and only if 
it is either $\delta_-$-stable or there exists a unique filtration of the form 
\be\label{eq:minusfiltr} 
0\subset \CE'\subset \CE 
\ee
so that $\CE'$ 
a semistable Higgs sheaf of rank $r(\CE')\geq 1$, and $\CE''=\CE/\CE'$ is a 
$\delta_-$-stable ADHM sheaf or rank $ r(\CE'')\geq 1$satisfying
\be\label{eq:hnslopecondB}
\mu(\CE') > \mu_{\delta_-}(\CE'') \qquad 
\mu(\CE') = \mu_{\delta_i}(\CE'').
\ee
\end{lemm}

{\it Proof.} 
It suffices to prove statement $(i)$ since the proof of $(ii)$ is 
analogous. 

Let $\delta_+\in (\delta_i, \ \delta_{i+1})$ be an arbitrary 
 noncritical stability parameter 
of type $(r,e)$.  
Suppose $\CE$ is a $\delta_i$-semistable ADHM 
sheaf on $X$. Then $\CE$ is either 
$\delta_+$-stable or there is a Harder-Narasimhan filtration of $\CE$ with respect to 
$\delta_+$-stability
\be\label{eq:plusHNfiltr}
0\subset \CE_1 \subset \cdots \subset \CE_h = \CE 
\ee
where $h\geq 2$. It is straightforward to check that $\CE_l$, $1\leq l\leq h$ 
must have rank $r(\CE_l)\geq 1$ and the successive quotients 
$\CE_{l+1}/\CE_l$, $0\leq l \leq h-1$ must also have rank $r(\CE_{l+1}/\CE_l)\geq 1$.  
Then by the general properties of Harder-Narasimhan filtrations 
\be\label{eq:HNineqB}
\mu_{\delta_+}(\CE_1)> \mu_{\delta_+}(\CE_2/\CE_1)>\cdots >
\mu_{\delta_+}(\CE_h/\CE_{h-1})
\ee
and 
\be\label{eq:ineqAB}
\mu_{\delta_+}(\CE_l)> \mu_{\delta_+}(\CE)
\ee
for all $1\leq l\leq h-1$. Since $\CE$ is $\delta_i$-semistable by assumption, 
inequalities \eqref{eq:ineqAB} imply that $v(\CE_l)=1$ for all $l=1,\ldots, h$. 
Therefore all quotients $\CE_{l+1}/\CE_{l}$, $0\leq l\leq h-1$ are semistable 
Higgs sheaves on $X$. 
Moreover, using the $\delta_i$-semistability condition and 
inequalities \eqref{eq:ineqAB} 
we have 
\be\label{eq:ineqB}
\delta_+\left({1\over r}-{1\over r(E_l)}\right) 
< \mu(E_l) -\mu(E) \leq \delta_i\left({1\over r}-{1\over r(E_l)}\right)
\ee
for all $l=1,\ldots, h$. 

Now let $\gamma \in (\delta_i, \ \delta_{i+1})$ be a fixed stability 
parameter. Then we claim that the set of 
isomorphism classes of locally free ADHM sheaves $\CE'$ on $X$ satisfying 
condition 
\begin{itemize}
\item[$(\star)$] 
There exists a $\delta_i$-semistable 
ADHM sheaf $\CE$ of type $(r,e)$ and a stability parameter $\delta_+
\in (\delta_i,\ \gamma]$ so that $\CE'\simeq \CE_l$ for 
some $l\in \{0,\ldots, h\}$, where $0\subset \CE_1\subset\cdots\subset 
\CE_h=\CE$, 
$h\geq 1$, is the  Harder-Narasimhan filtration of $\CE$ 
with respect to $\delta_+$-stability. 
\end{itemize} 
 is bounded. 
In order to prove this claim note that 
 that under the current conditions, inequalities 
\eqref{eq:ineqB} imply 
\[
-\gamma< \mu(E_l)-\mu(E) < {\delta_i\over r}.
\]
Moreover the set of isomorphism classes of $\delta_i$-semistable 
ADHM sheaves of type $(r,e)$ is bounded according to lemma 
(\ref{boundedness}). Therefore the above claim follows from Grothendieck's 
lemma. 

Then it follows that the set 
of numerical types $(r',e')$ of 
locally free ADHM sheaves $\CE'$ satisfying  property $(\star)$ 
is finite. This implies that there exists $0< \epsilon_+<\gamma-\delta_i$ 
so that for any $\delta_+\in (\delta_i, \ \delta_i+\epsilon_+)$, and 
any $\delta_i$-semistable ADHM sheaf $\CE$ of type $(r,e)$
inequalities \eqref{eq:ineqB} can be satisfied only if 
\be\label{eq:equal}
\mu_{\delta_i}(\CE_l) = \mu_{\delta_i}(\CE)
\ee
for all $l=1,\ldots, h$. 
Hence also 
\[ 
\mu(\CE_l/\CE_{l-1}) = \mu_{\delta_i}(\CE)
\]
for all $l=2,\ldots, h$. Then \eqref{eq:HNineqB} implies 
that we must have $h=2$. 
Therefore for all $\delta_+\in (\delta_i, \ \delta_i+\epsilon_+)$, 
any locally free $\delta_i$-semistable ADHM sheaf $\CE$ of type $(r,e)$ is either 
$\delta_+$-stable or has a 
Harder-Narasimhan filtration with respect to $\delta_+$-stability 
of the form \eqref{eq:plusfiltr} so that $\CE'$, $\CE''=\CE/\CE'$ 
satisfy conditions \eqref{eq:hnslopecondA}. 

Next note that the set of numerical types 
\be\label{eq:slopeset} 
{\sf {Sat}}_{\delta_i}(r,e) = \{(r',e')\in \IZ_{\geq 1}\times \IZ\, |\, r'\leq r ,\
r(e'+\delta_i) = r'(e+\delta_i)\}
\ee
is finite. Therefore $0< \epsilon_+<\gamma-\delta_i$ above may be chosen so 
that there are no critical stability parameters of type $(r',e')$ in the interval 
$(\delta_i, \ \delta_i+\epsilon_+)$ for any
$(r',e')\in {\sf {Sat}}_{\delta_i}(r,e)$. 
This implies in particular that for all $\delta_+\in (\delta_i, \ \delta_i+\epsilon_+)$, 
and any locally free $\delta_i$-semistable $\delta_+$-unstable 
ADHM sheaf $\CE$ of type $(r,e)$, the first step $\CE'$ in the Harder-Narasimhan 
filtration of $\CE$ must be $\delta_+$-stable rather than $\delta_+$-semistable. 

Conversely, suppose $\CE$ is a locally free ADHM sheaf of type $(r,e)$ 
on $X$ which has a 
filtration of the form \eqref{eq:plusfiltr} with $\CE'$ $\delta_+$-stable 
and satisfying conditions 
\eqref{eq:hnslopecondA}, for some 
$\delta_+\in (\delta_i, \ \delta_i+\epsilon_+)$. 
By the above choice of $\epsilon_+$, there are no critical stability parameters 
of type $(r(E'), d(E'))$ in the interval $(\delta_i, \ \delta_i+\epsilon_+)$
since $(r(E'), d(E'))\in  {\sf {Sat}}_{\delta_i}(r,e)$. 
Then  lemma (\ref{specstab}) implies that $\CE$ is $\delta_i$-semistable. 
This further implies that $\CE$ is $\delta_i$-semistable since it is an extension 
of semistable objects of equal $\delta_i$-slope.

\hfill $\Box$

Next we consider the behavior of $\delta$-semistable ADHM sheaves as 
$\delta$ specializes to $0$. The following results hold 
by analogy with  lemmas  (\ref{specstab}) and (\ref{minusdestabobj}) .
Since the proofs are very similar, they will be omited. 

\begin{lemm}\label{plusminuslemm} 
Let $X$ be a smooth projective curve over a field $K$ 
over $\IC$. Let  $\delta\in \IR\setminus \{0\}$
be a noncritical stability parameter of type $(r,e)\in \IZ_{\geq 1}\times \IZ$ so that 
there are no critical stability parameters of type $(r,e)$ in the 
interval $(0,\ \delta)$ if $\delta>0$, respectively 
$(\delta, \ 0)$ if $\delta<0$. 
Then any $\delta$-stable ADHM sheaf $\CE$ of type $(r,e)$ on $X$ 
is $0$-semistable. 
\end{lemm}

Conversely, 
\begin{lemm}\label{zeroHN} 
Under the same conditions let $(r,e)\in \IZ_{\geq 1}\times \IZ$. 
Then there exist $\epsilon_+>0$, $\epsilon_-<0$ so that the following hold. 

$(i)$ For any stability parameter $0<\delta_+ < \epsilon_+$ an
ADHM sheaf $\CE=(E,\Phi_1,\Phi_2,\phi,\psi)$ of type $(r,e)$ with $\psi$ nontrivial 
is $0$-semistable if and only if it is either $\delta_+$-stable or there exists 
a unique filtration $0\subset \CE'\subset \CE$ where $\CE'$ is a $\delta_+$-stable 
ADHM sheaf of rank $r(\CE')\geq 1$ and $\CE''=\CE/\CE'$ a semistable 
Higgs sheaf of rank $r(\CE'')\geq 1$ satisfying.
\be\label{eq:originsloepcondA}
\mu_{\delta_+}(\CE') > \mu(\CE''), \qquad \mu_0(\CE') = \mu(\CE'') = \mu_0(\CE)
\ee
Moreover an ADHM sheaf $\CE=(E,\Phi_1,\Phi_2,\phi,0)$ of type $(r,e)$ 
is $0$-semistable if and only if there exists a unique filtration $0\subset \CE' \subset \CE$
where $\CE'\simeq O=(0,\IC,0,0,0,0)$ and $\CE''=\CE/\CE'$ is a semistable Higgs 
sheaf of type $(r,e)$.

$(ii)$ For any stability parameter $\epsilon_-<\delta_- < 0$ an
ADHM sheaf $\CE=(E,\Phi_1,\Phi_2,\phi,\psi)$ of type $(r,e)$ with $\phi$ nontrivial 
is $0$-semistable if and only if it is either $\delta_-$-stable or there exists 
a unique filtration $0\subset \CE'\subset \CE$ where $\CE'$ is a 
semistable Higgs sheaf of rank $r(\CE')\geq 1$ and $\CE''=\CE/\CE'$ a semistable 
ADHM sheaf of rank $r(\CE'')\geq 1$ satisfying.
\be\label{eq:originsloepcondAB}
\mu(\CE') > \mu_{\delta_-}(\CE''), \qquad \mu(\CE') = \mu_0(\CE'') = \mu_0(\CE)
\ee
Moreover an ADHM sheaf $\CE=(E,\Phi_1,\Phi_2,0,\psi)$ of type $(r,e)$ 
is $0$-semistable if and only if there exists a unique filtration $0\subset \CE' \subset \CE$
where $\CE'$ is a semistable Higgs 
sheaf of type $(r,e)$  and 
$\CE''=\CE/\CE'\simeq O=(0,\IC,0,0,0,0)$.
\end{lemm}

\section{Moduli stacks and torus actions}\label{adhmmod}

The main goal of this section is to prove theorems (\ref{modthm}) and 
(\ref{virsmooth}). Other moduli stacks needed in the 
second part of this paper will be constructed. as well.
Natural torus actions on these stacks will be defined and  some 
structure results for the fixed loci will be proven.  
In the following $X$ is a 
smooth projective curve over $\IC$ and 
$M_1,M_2$ are fixed line bundles on $X$ equipped 
with a fixed isomorphism $M_1\otimes_X M_2 \simeq K_X^{-1}$ 
as in section (\ref{adhmstab}) as well as $E_\infty =\CO_X$.
The triple $(X,M_1,M_2)$ will be denoted by $\CX$.

\subsection{Moduli spaces of $\delta$-semistable ADHM sheaves}\label{constrmod}
Let $\CS$ be the category of schemes of finite type 
over $\IC$. Let $(r,e)\in \IZ_{\geq 1}\times \IZ$ 
be a fixed type and $\delta\in \IR_{>0}$ be a 
fixed stability parameter. Standard 
arguments show that flat families of 
$\delta$-semistable ADHM sheaves parameterized by 
complex schemes of finite type form a groupoid $\mfm^{ss}_\delta(\CX,r,e)$
over $\CS$.  For completeness recall, \cite[Def. 3.1]{modADHM},
that a flat family of ADHM sheaves on $X$ 
parameterized by a scheme $S$ of finite type over $\IC$ is an ADHM sheaf 
$\CE_S=(E_S, \Phi_{S,1,2}, \phi_S,\psi_S)$ on $X_S$ with twisting data 
$(M_1)_S$, $(M_2)_S$ and framing data $E_{S,\infty} = \CO_{X_S}$
so that $E_S$ is flat over $S$. An isomorphism of flat families of ADHM sheaves 
parameterized by $S$ is an isomorphism of ADHM sheaves on $X_S$.

{\it Proof of Theorem (\ref{modthm}).}
 Lemma
 (\ref{boundedness}) implies by standard arguments that the 
 groupoid $\mfm^{ss}_\delta(\CX,r,e)$ is an algebraic stack of finite type over $\IC$. 
 In particular, employing standard constructions for 
decorated sheaves as in 
\cite{semistpairs,stpairs,framed,decorated,modquivers,tensors}, 
 the groupoid $\mfm^{ss}_\delta(\CX,r,e)$ is shown to be isomorphic to a quotient stack 
 of a quasi-projective parameter scheme by a reductive affine 
 group.  
 
 For noncritical $\delta$, Corollary (\ref{automADHM}) implies that 
 all closed points have trivial stabilizers since there are no strictly semistable objects. 
 Then the moduli stack $\mfm^{ss}_\delta(\CX,r,e)$ is an algebraic space. 
The proof that $\mfm^{ss}_\delta(\CX,r,e)$ is isomorphic to a quasi-projective scheme 
in this case will be based on Theorem (\ref{adhmquivmod}) and 
Lemma (\ref{adhmquivstab}). 

Note that in each stability 
chamber there exists a rational stability parameter $\delta \in \IQ_{>0}$. 
For any $\sigma \in \IQ_{>0}$ let ${\mathfrak Q}^{ss}_{(\sigma,\delta)}(\CX,r,e;,1,0)$,
$\CQ^{ss}_{(\sigma,\delta)}(\CX,r,e;,1,0)$ be the 
moduli stack, respectively space of $(\sigma,\delta)$-semistable 
ADHM quiver sheaves on $X$ of type $(r,e;1,0)$. As observed 
in Remark (\ref{GITconstr}), there is a natural morphism 
\[
{\mathfrak q}: {\mathfrak Q}^{ss}_{(\sigma,\delta)}(\CX,r,e;,1,0)\to 
\CQ^{ss}_{(\sigma,\delta)}(\CX,r,e;,1,0).
\]
For sufficiently generic  $\sigma$ there are no strictly $(\sigma,\delta)$-semistable 
objects, and the above morphism is a $\IC^\times$-gerbe. 

Next note that there is also a natural morphism of algebraic stacks 
\be\label{eq:detmorphism} 
{\mathfrak d}_\infty: {\mathfrak Q}^{s}_{(\sigma,\tau)}(\CX,r,e; 1,0)\to 
{\mathfrak{Pic}}_{0}(X)
\ee
defined by taking the determinant line bundle of $E_\infty$, 
where 
${\mathfrak{Pic}}_{0}(X)$ is the algebraic stack of degree $0$ 
line bundles on $X$. 
Let ${\mathfrak Q}^{ss}_{(\sigma,\delta)}(\CX,r,e; 1,0)_{(0)}$ denote 
the stack-theoretic fiber of the morphism \eqref{eq:detmorphism} over the closed 
point of ${\mathfrak {Pic}}_{0}(X)$ determined by the trivial line bundle $\CO_X$. 
Let  ${\mathcal Q}^{ss}_{(\sigma,\delta)}(\CX,r,e; 1,0)_{(0)}$ be its scheme theoretic image 
via the morphism ${\mathfrak q}$, which is a closed subscheme of 
${\mathcal Q}^{ss}_{(\sigma,\delta)}(\CX,r,e; 1,0)$. 
Moreover, 
${\mathfrak Q}^{ss}_{(\sigma,\delta)}(\CX,r,e; 1,0)_{(0)}$ is again  $\IC^\times$-gerbe over 
${\mathcal Q}^{ss}_{(\sigma,\delta)}(\CX,r,e; 1,0)_{(0)}$.

By construction, an object of ${\mathfrak Q}^{ss}_{(\sigma,\delta)}(\CX,r,e; 1,0)_{(0)}$
over a scheme $S$ is a flat family $(E_S,E_{S,\infty}, \Phi_{S,i}, \phi_S,\psi_S)$ 
of $(\sigma,\delta)$-stable ADHM quiver sheaves such that $E_{S,\infty}$ is an 
invertible sheaf on $X\times S$, isomorphic to $\CO_{X_s}$ along each fiber $X_s$, $s\in S$. 
An object of $\mfm^{ss}_\delta(\CX,r,e)$ over $S$ is an ADHM sheaf $\CE_S$, flat over $S$. 
Then Lemmas (\ref{equivstab}), (\ref{sigmataudelta}) imply  
 that for sufficiently small $\sigma$ there is an epimorphism of algebraic stacks 
\be\label{eq:stackproj} 
{\mathfrak Q}^{ss}_{(\sigma,\tau)}(\CX,r,e; 1,0)_{(0)}\to 
\mfm^{ss}_\delta(\CX,r,e) 
\ee
which maps a flat family $(E_S,E_{S,\infty}, \Phi_{S,i}, \phi_S,\psi_S)$ 
to the flat family of $\delta$-stable ADHM sheaves 
$$(E_S\otimes E_{S,\infty}^{-1}, \Phi_{S,i}\otimes 1_{E_{S,\infty}^{-1}}, 
\phi_{S}\otimes 1_{E_{S,\infty}^{-1}},\psi_S\otimes 1_{E_{S,\infty}^{-1}}).
$$
 It is straightforward to check that 
the morphism \eqref{eq:stackproj} is a $\IC^\times$-gerbe and that it admits a canonical 
section.  Then \cite[Lemm. 3.18]{Laumon}
implies that there is an isomorphism of algebraic spaces 
$\mfm^{ss}_\delta(\CX,r,e)\simeq  {\mathcal Q}^{ss}_{(\sigma,\tau)}(\CX,r,e; 1,0)_{(0)}$. 

   \hfill $\Box$

 \begin{rema}\label{polarization} 
 As observed in remark (\ref{GITconstr}), the quasi-projective 
 moduli scheme ${\mathcal Q}^{ss}_{(\sigma,\delta)}(\CX,r,e; 1,0)$ 
 is naturally equipped with a polarization provided by the GIT construction 
 \cite[Thm. 1.10]{GIT}, \cite[Rem. 1.4.3.9]{GIT-decorated}. More precisely the 
G-linearized line bundle $L$ on the parameter space $P$ determines an ample 
line bundle $\CL$ on the GIT quotient. Therefore the closed subscheme 
${\mathcal Q}^{s}_{(\sigma,\delta)}(\CX,r,e; 1,0)_{(0)}$ is also equipped with a 
an ample line bundle $\CL_0$ obtained by restriction. 
   \end{rema}
 
 \subsection{Other moduli stacks}\label{biggerstack}
 For future reference we next construct 
 several moduli stacks of objects of the abelian category $\CC_X$, 
 where $X$ is a smooth projective curve over $\IC$ as in the previous 
 subsection. 
 
 Let $\delta \in \IR\setminus \{0\}$ be a stability parameter.
  Again, using standard arguments we construct the following groupoids over $\CS$ 
 \begin{itemize}
 \item $\obj(\CX)$: the groupoid of flat families of locally free 
 objects of $\CC_X$. An object of ${\mathfrak {Ob}}(\CX)$ over a $\IC$-scheme 
 $S$ of finite type is a flat family $\CE_S=(E_S, E_{S,\infty}, \Phi_S, \phi_S, \psi_S)$ 
 of locally free ADHM quiver sheaves of $X_S$, where 
 $E_{S,\infty}= \pi_S^*F_S$
 for some locally free $\CO_S$-module $F_S$ and the restriction 
 $\CE_S|_{X_s}$ for any point $s\in S$ is a locally free object of $\CC_{X_s}$
 The definition of isomorphisms is standard.  
 \item[$\bullet$]  ${\mathfrak {Ob}}(\CX)_{\leq 1}$: the groupoid of all flat families 
 of locally free objects of $\CC_X$ with $0\leq v\leq 1$. 
 The construction is the same as above, 
 except that $v$ may take only values in  $\{0,1\}$. Therefore the $\CO_S$-module 
 $F_S$ in the previous definition is either the zero module or an invertible 
sheaf on $S$. 
\item[$\bullet$] ${\mathfrak {Ob}}(\CX)_{v}$, $\obj(\CX,r,e,v)$: the groupoid of flat families 
 of locally free objects of $\CC_X$ with fixed $v\in \IZ_{\geq 0}$, respectively the 
 groupoid of flat families of locally free objects of $\CC_X$ with fixed type 
 $(r,e,v)\in \IZ_{\geq 1}\times \IZ\times \IZ_{\geq 0}$. 
   \item[$\bullet$] 
 ${\mathfrak {Ob}}_\delta^{ss}(\CX,r,e,v)$: the groupoid of flat families of 
 $\delta$-semistable objects of $\CC_X$ of fixed type $(r,e,v)\in \IZ_{\geq 1}
 \times \IZ\times \{0,1\}$, with $\delta \in \IR$. 
 \item ${\mathfrak{Ex}}(\CX)$: the groupoid of three term exact sequences of locally 
 free objects of $\CC_X$ constructed by analogy with \cite[Def. 7.2]{J-I}. 
 Note that there are canonical forgetful morphisms 
${\mathfrak p}, {\mathfrak p}', {\mathfrak p}'': {\mathrm{Ex}}(\CX)\to 
\obj(\CX)$.
 \item ${\mathfrak{Ex}}(\CX,\alpha,\alpha',\alpha'')$: the groupoid of three term exact sequences of locally 
 free objects of $\CC_X$ with fixed types $\alpha=(r,e,v)$, $\alpha'=(r',e',v')$ 
 $\alpha''=(r'',e'',v'')$ in $\IZ_{\geq 1}\times \IZ\times \IZ_{\geq 0}$, 
 $\alpha=\alpha'+\alpha''$. 
  \end{itemize}
 
 \begin{rema}\label{stackrem} 
 Note that the stack $\obj(\CX)_0$ is canonically isomorphic to the 
 stack ${\mathfrak {Higgs}}(\CX)$ of all locally free Higgs sheaves on $X$
 as defined in (\ref{Higgs-sheaves}). The stack $\obj^{ss}_\delta(\CX,r,e,0)$ 
 is canonically isomorphic to the moduli stack ${\mathfrak {Higgs}}^{ss}(\CX,r,e)$
 of semistable Higgs sheaves of 
 type 
$(r,e)$ on $X$ for any value of $\delta \in \IR\setminus\{0\}$. Both notations will be 
used interchangeably from now on. 
\end{rema} 
 
 The main properties of the above moduli stacks are summarized in the following lemma.   
  \begin{lemm}\label{allstack}
 $(i)$ For any $\delta \in \IR$  and any $(r,e,v)\in \IZ_{\geq 1}
 \times \IZ \times \{0,1\}$, ${\mathfrak {Ob}}_\delta^{ss}(\CX,r,e,v)$  
 is an Artin stack of finite type over $\IC$. 
 
 $(ii)$ The groupoids $\obj(\CX)$, 
 $\obj(\CX)_{\leq 1}$, $\obj(\CX)_{v}$, $\obj(\CX,r,e,v)$, $r,e,v\in
  \IZ_{\geq 1} \times \IZ \times \IZ_{\geq 0}$ are 
 Artin stacks of locally finite type 
over $\IC$ and there are open and closed immersion of Artin stacks 
\be\label{eq:opimmX}
\obj(\CX,r,e,v) \hookrightarrow 
\obj(\CX)_{v}\hookrightarrow \obj(\CX)_{\leq 1}\hookrightarrow \obj(\CX)\ee
for any $r,e,v\in
  \IZ_{\geq 1} \times \IZ \times \IZ_{\geq 0}$.
  
  $(iii)$
For any stability parameter $\delta\in \IR$ and 
any $(r,e)\in \IZ_{\geq 1} \times \IZ$ there are 
open immersions of Artin stacks 
\be\label{eq:opimmA} 
\begin{aligned}
& \obj^{ss}_\delta(\CX, r,e,1) \hookrightarrow \obj(\CX)_{1}\hookrightarrow 
\obj(\CX)_{\leq 1}\\
&  \obj^{ss}_\delta(\CX, r,e,0)\hookrightarrow \obj(\CX)_{0} \hookrightarrow 
\obj(\CX)_{\leq 1}.
\end{aligned}
\ee

$(iv)$ Suppose $\delta \in \IR$ is a noncritical stability 
parameter of type $(r,e)\in \IZ_{\geq 1} \times \IZ$. Then there is a morphism 
\be\label{eq:quotstackA} 
\obj_\delta^{ss}(\CX,r,e,1)\simeq \mfm_{\delta}^ss(\CX,r,e)
\ee
which identifies $\obj_\delta^{ss}(\CX,r,e,1)$ with a trivial 
$\IC^\times$-gerbe over $\mfm_{\delta}^ss(\CX,r,e)$. 

$(v)$ The groupoids ${\mathfrak{Ex}}(\CX)$, 
${\mathfrak{Ex}}(\CX,\alpha,\alpha',\alpha'')$, $\alpha,\alpha',\alpha''
\in \IZ_{\geq 1}\times \IZ\times \IZ_{\geq 0}$ are algebraic stacks of locally 
finite type over $\IC$ and the canonical forgetful morphisms 
${\mathfrak p}, {\mathfrak p}', {\mathfrak p}'': {\mathfrak{Ex}}(\CX)\to 
\obj(\CX)$ are of finite type over $\IC$. 
\end{lemm}

{\it Proof.}
The proof of lemma (\ref{allstack}.$i$) follows from the boundedness 
lemma (\ref{boundedness}), again in complete analogy with 
\cite[Thm 1.1]{modADHM}.

Given the construction of the parameter spaces for Higgs sheaves 
\cite{semistpairs,simpsonI,projective} and for ADHM sheaves \cite{modADHM}, 
statements (\ref{allstack}.$ii$) and (\ref{allstack}.$v$) 
follow by analogy with \cite[Thm 4.6.2.1]{Laumon}, 
\cite[Thm 9.4]{J-I}, \cite[Thm. 9.6]{J-I}. 

The fact that the natural morphisms \eqref{eq:opimmA} are open 
immersions follows from the fact that Higgs sheaf semistability as well 
as $\delta$-semistability are open conditions in flat families. 

Lemma (\ref{allstack}.$iv$) follows from lemmas (\ref{equivstab}) and (\ref{endlemma})
taking into account remark (\ref{isomrema}). 

\hfill $\Box$

\subsection{Torus actions and fixed foci}\label{toractsect} 

Next we define natural torus actions on the above moduli stacks. 
We will employ the definition of group actions on   stacks 
given in \cite[Def.1.3, Def.2.1]{gract-stacks}. 
Let 
${\bf T}=\IC^\times\times \IC^\times$. Then 
  there is a torus action 
 ${\bf T}^\times \times \obj(\CX) \to \obj(\CX)$ 
 defined as follows. Given any flat family $\CE_S=(E_S,(\pi_S^*F_S)^{\oplus v},
 \Phi_{S,1},\Phi_{S,2},
\phi_S,\psi_S)$ parameterized by a scheme $S$ of finite type 
over $\IC$, and $t_1,t_2:S\to {\bf T}$, set 
\[
(t_1,t_2)\times \CE_S \to \CE_S^{(t_1,t_2)} 
\]
where 
\be\label{eq:Tact}
\CE_S^{(t_1,t_2)}= 
(E_S,(\pi_S^*L_S)^{\oplus v},t_1\Phi_{S,1},t_2\Phi_{S,2},t_1t_2\phi_S,\psi_S).
\ee
Moreover, if $\xi : \CE_S{\buildrel \sim \over \longto}  \CE'_S$ is an isomorphism 
of flat families over $S$, the isomorphism $\xi^{(t_1,t_2)} 
: \CE_S^{(t_1,t_2)}{\buildrel \sim \over \longto}  {\CE'}_S^{(t_1,t_2)}$
is given by the same isomorphism $\xi:E_S \to E_S$ of coherent 
$\CO_{X_S}$-modules, since ${\bf T}$
acts linearly on the data $(\Phi_{S,1,2},\phi_S,\psi_S)$ leaving the underlying coherent 
sheaf $E_S$ unchanged. 

Let ${\bf S}\simeq \IC^{\times}\subset {\bf T}$ be the antidiagonal subtorus defined 
by the embedding $t\to (t^{-1},t)$. Then the action \eqref{eq:Tact} 
induces an {\bf S}-action on $\obj(\CX)_{\leq 1}$. In this case we will use the notation  $
\CE_S^t=\CE_S^{(t^{-1},t)}$, $\xi^t = \xi^{(t^{-1},t)}$. 

Obviously, there are analogous actions 
on any  substack of the form $\obj(\CX)_{\leq 1}$, $\obj(\CX)_v$, $\obj(\CX,r,e,v)$
 $\obj_\delta^{ss}(\CX,r,e,v)$, with $(r,e,v)\in\IZ_{\geq 1}\times \IZ\times \IZ_{\geq 0}$
 so that 
the open immersions \eqref{eq:opimmX}, \eqref{eq:opimmA} are equivariant.
Moreover, there are analogous torus actions ${\bf T}\times \mfm_\delta^{ss}
(\CX,r,e) \to \mfm_\delta^{ss}(\CX,r,e)$, ${\bf S}\times \mfm_\delta^{ss}(\CX,r,e) \to 
\mfm_\delta^{ss}(\CX,r,e)$ obtained by setting $v=1$ and $F_S=\CO_S$ in 
\eqref{eq:Tact}. 

Finally note that since the {\bf T} and {\bf S} actions defined above are linear 
on the ADHM data, it is straightforward to check that they can be canonically 
lifted to 
torus actions on the stack ${\mathfrak {Ex}}(\CX)$ so that the forgetful 
morphisms ${\mathfrak p},{\mathfrak p}',{\mathfrak p}''$ are equivariant. 
A flat family 
\[
0\to \CE_S'{\buildrel \xi_S\over \longto} \CE_S {\buildrel \eta_S\over 
\longto} \CE_S''\to 0 
\]
of three term exact sequences is mapped by $(t_1,t_2):S \to {\bf T}$ 
to 
\[
\xymatrix{
0\ar[r]&  \CE_S'^{(t_1,t_2)}\ar[r]^-{\xi_S^{(t_1,t_2)}}& \CE_S^{(t_1,t_2)} 
\ar[r]^-{\eta_S^{(t_1,t_2)}}& 
\CE_S''^{(t_1,t_2)}\ar[r]& 0 \\}
\]
where $\xi_S^{(t_1,t_2)}$, $\eta_S^{(t_1,t_2)}$ have the same underling 
morphisms of $\CO_{X_S}$-modules as $\xi_S,\eta_S$. Obviously, this torus 
action preserves the substacks of fixed type $(\alpha,\alpha',\alpha'')$. 

Let $\mfm$ denote one of the stacks  in lemma (\ref{allstack}.$i$) - (\ref{allstack}.$iv$).
 Let $\mfm^{\bf S}$ be the stack theoretic fixed locus as defined in 
\cite[Prop 2.5]{gract-stacks}. Recall that a morphism $S\to \mfm^{\bf S}$, 
where $S$ is a scheme of finite type over $\IC$, is determined by the data 
$\{\CE_S, \xi_S(t)\,|\, t:S\to {\bf S}\}$ where  
$\CE_S$ is a flat family of locally free objects of $\CC_X$ parameterized by $S$ 
and for any morphism $t:S\to {\bf S}$
\be\label{eq:fixedisomA}
\xi_S(t): \CE_S{\buildrel \sim\over \longto} \CE_S^t
\ee
is an isomorphism of flat families over $S$ satisfying the identity 
\be\label{eq:fixedisomB}
\xi_S({t't}) = \xi_S({t})^{t'} \circ \xi_S(t') 
\ee
for all $t,t': S \to {\bf S}$ (see  \cite[Prop 2.5]{gract-stacks}). 
In the following a flat family satisfying this property will be called 
{\bf S}-fixed up to isomorphism. 

\begin{lemm}\label{fixedlociA} 
Let $\CE_S=(E_S, \pi_S^*F_S, \Phi_{S,1,2}, \phi_S, \psi_S)$ 
be a ${\bf S}$-fixed flat family of locally free 
objects of $\CC_X$ parameterized by a connected scheme 
$S$ of finite type over $\IC$. 
Then there are direct sum decompositions 
\be\label{eq:fixedcondA} 
E_S \simeq \oplus_{n\in I} E_S(n), \qquad F_S = \oplus_{n\in J} F_S(n)
\ee
with $ I,J\subset \IZ$ finite subsets,  satisfying the following conditions 
\begin{itemize} 
\item If $E_S$, respectively $F_S$, are not the zero sheaf, no direct summand in 
\eqref{eq:fixedcondA} is the zero sheaf.
\item All components $\Phi_{S,i}(n,m) : E_S(n)\otimes_{X_S}(M_i)_S \to 
E_S(m)$, $n,m\in I$, with respect to the direct sum decomposition 
\eqref{eq:fixedcondA} are identically zero if $m\neq n + (-1)^{i-1}$ for $i=1,2$. 
\item All components $\phi_S(n,m) : E_S(n)\otimes_{X_S}M_S \to \pi_S^*F_S(m)$, 
$n\in I$, $m\in J$, with respect to the direct sum decomposition 
\eqref{eq:fixedcondA} are identically zero if $m\neq n$. 
\item All components $\psi_S(n,m) : \pi_S^*F_S(n) \to E_S(m)$, 
$n\in I$, $m\in J$, with respect to the direct sum decomposition 
\eqref{eq:fixedcondA} are identically zero if $m\neq n$. 
\end{itemize}
\end{lemm}

{\it Proof.}  Note that $\xi_S(t)^{t'}$ 
in equation   \eqref{eq:fixedisomB} is identical to $\xi_S(t)$ 
as an isomorphism of $\CO_{X_S}$-modules, 
as observed below \eqref{eq:Tact}. This implies that the underlying 
coherent $\CO_{X_S}$-module $E_S$ of a flat family ${\CE}_S$ 
fixed by ${\bf S}$ up to isomorphism has a ${\bf S}$-linearization. Then 
lemma (\ref{fixedlociA})  follows from a standard analysis of the fixed locus conditions
analogous to \cite[Prop. 3.15]{modADHM}. Details will be omitted. 

\hfill $\Box$ 

A final technical result need below is the following. 
Recall that in the proof of Theorem (\ref{modthm}) the moduli 
space $\mfm_\delta^{ss}(\CX,r,e)$, for noncritical $\delta$, 
has been identified with the quasi-projective scheme 
${\mathcal Q}_{(\sigma,\delta)}^{ss}(\CX,r,e;,1,0)_{(0)}$. 
According to Remark (\ref{polarization}), the latter is equipped with an ample 
line bundle $\CL_0$ determined by the GIT construction of the ambient 
moduli scheme ${\mathcal Q}_{(\sigma,\delta)}^{ss}(\CX,r,e;,1,0)$.
Given the construction of the GIT parameter space $P$ in \cite{modquivers}, it is easy to check 
that the torus action 
defined in \eqref{eq:Tact} is induced by a natural torus action on $P$ which commutes with the action of the GIT group $G$. Moreover, the torus action on $P$ admits a natural lift to a linearization 
on the GIT line bundle $L$. This implies that the line bundle $\CL_0$ is {\bf T}-equivariant. 
Then one can use the same  proof as in  \cite[Cor 2]{equiv-compl} to conclude: 
\begin{lemm}\label{equivcover} 
There exists a Zariski open cover of $\mfm_\delta^{ss}(\CX,r,e)$ consisting of finitely many 
{\bf T}-invariant affine charts. 
\end{lemm}

 \subsection{Virtual smoothness for noncritical $\delta$}\label{torvirt}
 In this subsection we prove theorem (\ref{virsmooth}). By analogy with 
the proof of theorem \cite[Thm. 1.5]{modADHM} we first need a 
vanishing result in the deformation theory of $\delta$-stable 
ADHM sheaves. 
According to \cite[Prop 4.5, 4.9]{modADHM} the deformation complex
$\CC(\CE)$ of 
a locally free ADHM sheaf $\CE$ on $X$ is obtained from the complex 
$\CC(\CE,\CE)$ defined in Proposition (\ref{quivextA}), equation \eqref{eq:hypercohA}
by removing the direct summand ${{Hom}}_X(E_\infty,E_\infty)$ from the first 
term and setting $\alpha_\infty=0$ in the expression of $d_1$. The resulting complex, 
$\CC_{\mathrm{def}}(\CE)$ will be referred to as the deformation complex of $\CE$. 

Then the main technical element in the proof of virtual smoothness 
is the following vanishing result for the hypercohomology groups 
of the resulting complex $\CC_{\mathrm{def}}(\CE)$.  The proof is similar to the proof of lemma 
\cite[Lemma 4.10]{modADHM}. Details will be omitted. 
\begin{lemm}\label{hypervanA}
Let $\CE$ be a $\delta$-stable ADHM sheaf on $X$. Then 
$\IH^i(X,\CC_{\mathrm{def}}(\CE))=0$, for all $i\geq 3$ and for all $i\leq 0$. 
\end{lemm} 

{\it Proof.} 
The proof is similar to the proof of lemma 
\cite[Lemma 4.10]{modADHM}. Using Serre duality it reduces to proving 
that the automorphism group of a $\delta$-stable ADHM sheaf is trivial. 
Details will be omitted. 

\hfill $\Box$. 

Let 
${\bf T}=\IC^\times\times \IC^\times$ and ${\bf S}\subset {\bf T}$ act on 
$\mfm_\delta^{ss}(\CX,r,e)$ as in \eqref{eq:Tact}. 
These are the  torus actions employed in theorem (\ref{virsmooth}) whose proof will be 
presented below. 

{\it Proof of theorem (\ref{virsmooth}).} 
If $\delta \in \IR_{>0}$, is a noncritical stability parameter of type 
$(r,e)\in \IZ_{\geq 1}\times \IZ$, any $\delta$-semistable ADHM sheaf 
of type $(r,e)$ on $X$ is $\delta$-stable. 
Given lemmas (\ref{equivcover}),  (\ref{hypervanA}), the existence of a 
{\bf T}-equivariant 
perfect tangent-obstruction theory follows by repeating the steps 
in \cite[Sect 5.]{modADHM} in the present setting. Details will 
be omitted. Obviously, the resulting perfect tangent-obstruction theory 
is also ${\bf S}$-equivariant.

In order to prove the second part of theorem (\ref{virsmooth}), 
note that there is a universal ADHM locally free ADHM sheaf on 
$\mfm_\delta^{ss}(\CX,r,e)\times X$ since all stable objects have trivial 
automorphisms. Let ${\mathfrak p}: \mfm_\delta^{ss}(\CX,r,e)\times X
\to \mfm_\delta^{ss}(\CX,r,e)$, $\pi_X:  \mfm_\delta^{ss}(\CX,r,e)\times X
\to X$ be the canonical projections. Let $\CC_{\mathrm{def}}({\mathfrak E})$ be the 
deformation complex of the universal object ${\mathfrak E}$. 
Then 
Grothendieck duality \cite{nironi} for the projection morphism ${\mathfrak p}$ 
yields an isomorphism 
\[
({\bf R}{\mathfrak p}_*(\CC({\mathfrak E})))^\vee 
\simeq {\bf R}{\mathfrak p}_* {\bf R}\lochom(\CC({\mathfrak E} ), 
\pi_X^*K_X[1]). 
\]
This isomorphism is compatible with the induced {\bf T} as well
as ${\bf S}$ actions.

Taking into account the isomorphism $M\simeq K_X^{-1}$ 
and the fact that $\CC({\mathfrak E})$ is locally free, 
a straightforward computation shows that 
\[
{\bf R}\lochom(\CC({\mathfrak E}), \pi_X^*K_X[1])\simeq 
\CC({\mathfrak E})[1]
\]
as ${\bf S}$-equivariant complexes. 
Therefore we obtain an isomorphism of ${\bf S}$-equivariant 
complexes 
\[
({\bf R}{\mathfrak p}_*(\CC({\mathfrak E})))^\vee 
\simeq {\bf R}{\mathfrak p}_*(\CC({\mathfrak E}))[1].
\]
This yields the required 
${\bf S}$-equivariant 
nondegenerate pairing  on the perfect tangent-obstruction theory of 
$\mfm_\delta(\CX,r,e)$ \cite{micro,symm}.

 \hfill $\Box$ 
 
\section{Properness of torus fixed loci}\label{fixedsect}
The main goal of this section is to prove theorem 
(\ref{Sfixed}) working under the same conditions 
as in section (\ref{adhmmod}).
Recall that in the proof of Theorem (\ref{modthm}) the moduli space
$\mfm_\delta^{ss}(\CX,r,e)$, with noncritical $\delta$, has been 
identified with the  closed subscheme ${\mathcal Q}^{ss}_{(\delta,\sigma)}(\CX,r,e;1,0)_{(0)}$
of a quasi-projective 
moduli scheme ${\mathcal Q}^{ss}_{(\delta,\sigma)}(\CX,r,e;1,0)$ 
as in  Theorem (\ref{adhmquivmod}) for sufficiently small $\delta$. 
Recall also that there is a proper Hitchin morphism 
\be\label{eq:adhmhitchinmap} 
h_Q : {\mathcal Q}^{ss}_{(\delta,\sigma)}(\CX,r,e;1,0)\to \IV 
\ee
where $\IV$ is an affine space. 
Since the fixed locus 
$\mfm_\delta^{ss}(\CX,r,e)^{\bf S}$ is a closed subscheme of 
$\mfm_\delta^{ss}(\CX,r,e)^{\bf S}$
in order to prove that it is proper over $\IC$ it suffices to prove that it is  a subscheme of the 
fiber of the Hitchin morphism at $0\in \IV$. In order to establish this it suffices to prove that 
the polynomial invariants of $\IC$-valued ${\bf S}$-fixed points vanish since vanishing 
at closed points implies vanishing in families. 
The proof will proceed on step-by-step basis, as shown below. 

\subsection{{\bf S}-fixed semistable Higgs sheaves}\label{fixedhiggssect}
One first  proves the analogous result for semistable Higgs sheaves $\CE=(E,\Phi_{1},\Phi_2)$ on $X$. 
More precisely, recall that there is a proper Hitchin map 
$h: Higgs^{ss}(\CX,r,e) \to \IH$ where 
 $$\IH = \oplus_{n=0}^r H^0(X, \mathrm{Symm}^n(M_1^{-1}\oplus M_2^{-1}) ),$$
 mapping a polystable Higgs sheaf of type $(r,e)$ to its characteristic polynomial. 
 In order to simplify the exposition, we will adopt the following notation conventions
for Higgs sheaves. Suppose $\CE=(E,\Phi_1,\Phi_2)$ is a Higgs sheaf on $X$ 
with coefficient sheaves $M_1,M_2$. 
Given a finite sequence $\Phi_{i_n}: E\otimes_X M_{i_n} \to E$, where 
$n\geq 1$ and $(i_1,\ldots, i_n) \in \{1,2\}^{\times n}$, the composition 
\[
\Phi_{i_n} \circ (\Phi_{i_{n-1}} 
\otimes 1_{M_{i_n}}) \circ \cdots \circ (\Phi_{i_1} \otimes 1_{M_{i_2}\otimes_X
\cdots \otimes_X M_{i_n}}): E\otimes_X M_{i_1}\otimes_X \ldots 
\otimes_X M_{i_n} \to E 
\]
will be denoted by 
\[
\Phi_{i_n} \Phi_{i_{n-1}}\ldots \Phi_{i_1} : E\otimes_X M_{i_1}\otimes_X \ldots 
\otimes_X M_{i_n} \to E .
\]
\begin{lemm}\label{higgsfixed}
Let $\CE$ be an {\bf S}-fixed semistable Higgs sheaf of type $(r,e)$ on $X$. 
Then $h([\CE])=0$. 
\end{lemm}

 {\it Proof.}
 Let $\CE=(E,\Phi_{1},\Phi_{2})$.
 Given the construction of the Hitchin map, 
 the  image $h([\CE])$ is determined by the 
 polynomial invariants 
\be\label{eq:hitmap}
(\mathrm{tr}(\Phi_{1}^{n_1}\Phi_{2}^{n_2})), \qquad n_1,n_2\leq 0 , \quad
n_1+n_2\leq r.
\ee
Any {\bf S}-fixed semistable Higgs sheaf 
must satisfy conditions $(i)$ and $(ii)$ of lemma (\ref{fixedlociA}). 

Note that $\mathrm{tr} (\Phi_{1}^{n_1}\Phi_{2}^{n_2})$ 
is a homogeneous element of $\IH$ for the action of ${\bf S}$ with 
weight $n_1-n_2$. Therefore if $(E,\Phi_{1},\Phi_{2})$ is fixed by the 
${\bf S}$-action, 
we must have 
\be\label{eq:trvan}
\mathrm{tr} (\Phi_{1}^{n_1}\Phi_{2}^{n_2})=0 
\ee
for all $n_1,n_2\geq 0$, $0<n_1+n_2\leq r$, $n_1\neq n_2$. 
In order to prove the claim we have to show that the vanishing result 
\eqref{eq:trvan} holds for $n_1=n_2$ as well. 

Note that only finitely many terms in character decomposition of $E$ 
are nontrivial. Hence 
\be\label{eq:chardecomp}
E =\bigoplus_{s_1\leq n\leq s_2} E_{}(n)
\ee
for some $s_1,s_2\in \IZ$, $s_1\leq s_2$. If $s_1=s_2$, $\Phi_{1}$, 
$\Phi_{2}$ are trivial according to lemma (\ref{fixedlociA}) 
and there is nothing to prove. Therefore we will 
assume $s_1< s_2$. Then lemma (\ref{fixedlociA}) implies that only 
components of the form
\[
\Phi_{1}(n): E_{}(n)\to E_{}(n+1)\qquad \Phi_{2}(n-1): E(n)\to E(n-1)
\]
are allowed to be nontrivial. 
This implies that the monomials $(\Phi_{1}\Phi_{2})^n$, $n\geq 1$ have the following 
block form with respect to the decomposition \eqref{eq:chardecomp}
\[
({\Phi_{1}}{\Phi_{2}})^n = \mathrm{diag}\left(0, ({\Phi_{1}}(s_1){\Phi_{2}}
(s_1))^n, 
\ldots, ({\Phi_{1}}(s_2-1){\Phi_{2}}(s_2-1))^n\right).
\]
Using the structure results proven in lemma 
(\ref{fixedlociA}),  condition \eqref{eq:higgsrel} in 
definition (\ref{Higgs-sheaves}) is equivalent to the following relations 
\be\label{eq:higgsrelA}
\Phi_{1}(n)\Phi_{2}(n) = \Phi_{2}(n+1)\Phi_{1}(n+1)
\ee
for all $s_1\leq n\leq s_2-1$. In particular, 
\be\label{eq:higgsrelB}
\Phi_{1}(s_2-1)\Phi_{2}(s_2-1) =0,\qquad 
\Phi_{2}(s_1)\Phi_{1}(s_1)=0.
\ee
If $s_2=s_1+1$, the required vanishing result follows immediately from 
relations \eqref{eq:higgsrelB}. Hence we will assume $s_1\leq s_2-2$ 
in the following. 
Then we will prove by an inductive argument that 
\be\label{eq:higgsrelC}
\mathrm{tr}({\Phi_{1}}(s_2-k) {\Phi_{2}}(s_2-k))^n =0 
\ee
for all $1\leq k \leq (s_2-s_1)$ and all $n\geq 1$.
The case $k=1$ follows immediately from relations 
\eqref{eq:higgsrelB}. 
Then note that for any $2\leq k \leq (s_2-s_1)$ and any $n\geq 1$
we have 
\[
\begin{aligned} 
\mathrm{tr}({\Phi_{1}}(s_2-k){\Phi_{2}}(s_2-k))^n 
& = \mathrm{tr}({\Phi_{2}}(s_2-k+1){\Phi_{1}}(s_2-k+1))^n \\
& = \mathrm{tr}({\Phi_{1}}(s_2-k+1){\Phi_{2}}(s_2-k+1))^n\\
\end{aligned}
\]
using invariance of the trace under cyclic permutations of the arguments. 
This proves the inductive step, hence the required vanishing result follows.

\hfill $\Box$

   \subsection{Asymptotic fixed loci}\label{asympfixedsect}
Next we prove that a similar result holds for asymptotically stable
ADHM sheaves of type $(r,e)\in \IZ_{\geq 1}\times \IZ$ on $X$. 

\begin{lemm}\label{asympfixed}
Suppose $\delta>\delta_N$ is an asymptotic stability 
parameter. Then any $h([\CE])=0$ for any $\delta$-stable ADHM sheaf on $X$. 
In particular the {\bf S}-fixed locus 
$\mfm_\delta^{ss}(\CX,r,e)^{\bf S}$ is proper over  $\IC$.
\end{lemm}

{\it Proof.}
According to \cite[Lemma 2.5]{modADHM} any asymptotically stable 
ADHM sheaf on $X$ with $E_\infty =\CO_X$ has $\phi=0$. 
Since $\phi=0$, the Hitchin map maps the isomorphism class of an 
asymptotically stable ADHM sheaf to the characteristic polynomial 
of the underlying Higgs sheaf. 
Then the proof of Lemma (\ref{asympfixed}) is identical to the proof 
of lemma (\ref{higgsfixed}). 

\hfill $\Box$

The proof of Theorem (\ref{Sfixed}) can then be concluded by an inductive argument.

\subsection{Rank and chamber induction} 
The proof will proceed by induction on the rank $r\geq 1$. 

Lemmas (\ref{rankone}), (\ref{asympfixed}) imply that any rank one stable ADHM sheaf has 
trivial polynomial invariants. 

Let $r\geq 2$. 
Suppose any {\bf S}-fixed $\delta$-stable ADHM sheaf of type $(r',e)$ with $1\leq r'<r$, $e\in \IZ$ and, 
for a fixed type $(r',e)$, any noncritical 
stability parameter $\delta\in \IR_{>0}$ of type $(r',e)$ has trivial polynomial invariants. 
Then we have to prove the
analogous statement holds for each fixed type $(r,e)$, $e\in \IZ$. 
This will be done by a chamber inductive argument starting with the 
asymptotic chamber $\delta >\delta_N$. Lemma (\ref{asympfixed})
 implies that 
the statement holds 
for $\delta> \delta_N$. If $N=0$, there is nothing to prove, so we will 
assume $N\geq 1$ in the following. Set $\delta_0=0$ and 
$\delta_{N+1}= +\infty$. 

Suppose the statement holds for all {\bf S}-fixed $\delta$-stable ADHM sheaves 
of type $(r,e)$, where $\delta \in (\delta_{i},\ \delta_{i+1})$, for some $i=1,\ldots, N$. 
Then we will prove that the same holds for {\bf S}-fixed $\delta$-stable ADHM sheaves of type 
$(r,e)$ with $\delta \in (\delta_{i-1},\ \delta_i)$. 

Let $\CE$ be an {\bf S}-fixed $\delta_-$-stable ADHM sheaf of type $(r,e)$ on $X$
 for some $\delta_- \in (\delta_{i-1},\ \delta_i)$.
Then lemma (\ref{specstab}) implies that one of the 
following cases must hold 
\begin{itemize}
\item[$(1)$] 
$\CE$ is $\delta_i$-stable, hence also $\delta_+$-stable
for any $\delta_+\in (\delta_i, \ \delta_{i+1})$, or 
\item[$(2)$] $\CE$ is strictly $\delta_i$-semistable 
and there is a nontrivial extension 
\be\label{eq:casetwoext}
0\to \CE' \to \CE \to \CE''\to 0
\ee
in the abelian category $\CC_{X}$,
 where $\CE'$ is a $\delta_i$-stable 
ADHM sheaf on $X_K$, $\CE''$ is a semistable Higgs sheaf on $X$
with $\mu_{\delta_i}(\CE') = \mu(\CE'')$. 
\end{itemize}
In the first case, there is nothing to prove, hence suppose $(2)$ holds. 
Then Corollary (\ref{wallHN}) implies that $0\subset \CE'\subset \CE$ is a 
Harder-Narasimhan filtration with respect to $\delta_+$-stability for any 
$\delta_+\in (\delta_i, \ \delta_{i+1})$. 
 Since $\CE$ is {\bf S}-fixed, for each $t\in {\bf S}$ there is an isomorphism 
 $\xi({t}): \CE \to \CE^{t}$ as in equation \eqref{eq:fixedisomA}. 
 Moreover it straightforward to check that the 
 torus action \eqref{eq:Tact} preserves the Harder-Narasimhan filtrations i.e. 
 \[
 (\CE')^{t} = (\CE_{}^{t})'.
 \] 
 This implies that both $\CE',\CE''$ must be {\bf S}-fixed. 

Then the polynomial invariants of $\CE'$ are 
trivial according to the rank induction hypothesis and those 
of $\CE''$ are trivial according to Lemma (\ref{higgsfixed}). 
Since $\CE$ is an extension of quiver sheaves, this implies that its 
polynomial invariants are also trivial. 

\hfill $\Box$

\section{Versal deformations and holomorphic Chern-Simons theory}\label{versCS}

The main goal of this section is to show that the recent results 
\cite{genDTI,genDTII} of Joyce and Song on Behrend functions 
for algebraic moduli stacks of coherent sheaves on smooth projective 
Calabi-Yau threefolds also hold for the moduli stacks constructed in section
(\ref{biggerstack}). This will be needed in the second part of this paper, 
which is concerned with  wallcrossing formulas 
for ADHM invariants.

Let $\obj(\CX)_{\leq 1}$, $\obj(\CX)_{0}$ be 
the algebraic moduli stacks constructed in lemma (\ref{allstack}). 
Let $\CM(\CX)_{\leq 1}$, $\CM(\CX)_{0}$ denote the coarse algebraic moduli space of 
simple objects
in $\obj(\CX)_{\leq 1}$, $\obj(\CX)_{0}$. 
The first goal of this section is to prove that the statements of \cite[Thm. 5.2]{genDTI} 
and \cite[Thm. 5.3]{genDTI} hold for  
$\CM(\CX)_{\leq 1}$, $\obj(\CX)_{\leq 1}$, respectively $\obj(\CX)_{0}$, 
$\CM(\CX)_{0}$. 
For simplicity, we will denote 
extensions groups in the abelian category $\CC_X$ by $\rm{Ext}(\ ,\ )$ 
in the following. 
Then 
by analogy with \cite{genDTI} we claim 
\begin{theo}\label{theoremA} 
Let $\CM$ denote either $\CM(\CX)_{\leq 1}$ or $\CM(\CX)_{0}$.
Then for each $\IC$-valued point $[\CE]\in \CM(\IC)$ there exists a finite-dimensional 
complex manifold $U$, a holomorphic function $\omega:U\to \IC$ and a point $u\in U$ 
so that $\omega(u)=d\omega(u)=0$ and $\CM(\IC)$ is locally isomorphic as a complex 
analytic space to $\mathrm{Crt}(\omega)$ near $u$. Moreover, $U$ can be taken to be 
isomorphic 
to an open neighborhood of $u=0$ in the vector space $\mathrm{Ext}^1(\CE,\CE)$. 
\end{theo} 

Similarly, let $\mfm$ denote either $\obj(\CX)_{\leq 1}$or 
 $\obj(\CX)_{0}$. According to \cite[Thm 5.3]{genDTI}, the general theory of Artin stacks 
implies that 
for any $\IC$-valued point $[\CE]\in \mfm(\IC)$, there exists an $\mathrm{Aut}(\CE)$-
invariant subscheme $S\subset \mathrm{Ext}^1(\CE,\CE)$ 
over $\IC$ parameterizing an $\mathrm{Aut}(\CE)$-equivariant 
 versal family $\CE_S$ of locally free objects of $\CC_X$ with $v\in \{0,1\}$
 so that $\CE_S|_{X_0} \simeq \CE$. Moreover there exists an 
\'etale morphism of Artin stacks $\Phi: [S/\mathrm{Aut}(E)]\to 
\mfm$ so that $\Phi([0])= [\CE]$, 
\be\label{eq:natmorphA} 
\Phi_*: \mathrm{Stab}([0])\simeq 
\mathrm{Aut}(\CE) \to 
\mathrm{Stab}([\CE])
\ee
and 
\be\label{eq:natmorphB} 
d\Phi: T_{[0]}[S /\mathrm{Aut}(\CE)]\simeq \mathrm{Ext}^1(\CE,\CE)
\to T_{[\CE]}\mfm
\ee
are natural isomorphisms. Then the following holds. 

\begin{theo}\label{theoremB} 
For each $\IC$-valued point $[\CE]\in \mfm(\IC)$ there exists an 
open neighborhood  of $0$ $ U\subset \mathrm{Ext}^1(\CE,\CE)$ 
in the analytic topology, 
a holomorphic function $\omega:U\to \IC$ 
so that $\omega(0)=d\omega(0)=0$ 
and an open neighborhood of $0$ $V\subset S_{\sf an}$ 
so that there is an isomorphism of complex analytic spaces $\Xi: \mathrm{Crt}(\omega)\to 
V$ satisfying $\Xi(0)=0$ and $d\Xi_0=\mathrm{Id}_{\mathrm{Ext}^1(\CE,\CE)}$. 
Moreover, if $G$ is a maximal compact subgroup of $\mathrm{Aut}(\CE)$, 
$U,w,V$ can be chosen to be $G^\IC$-equivariant. 
\end{theo}

{\it Proof of Theorems (\ref{theoremA}), (\ref{theoremB}).}
Since we have restricted ourselves to locally free objects,
theorems (\ref{theoremA}), (\ref{theoremB}) can be proven using gauge theoretic 
methods in complete analogy with \cite[Thm 5.2]{genDTI}, 
\cite[Thm 5.3]{genDTI}. One has to check that the main arguments in 
\cite{M-Kuranishi}, \cite{genDTI} carry over to the current decorated bundle moduli problem. 
A gauge theoretic approach to various decorated bundle moduli problems has been 
previously employed for example in \cite{pairsI,pairsII,augmented,OT-sw,BDW,equiv,
T-nonab,OT-four,KL-moduli,HKpairs,HKquivers,triples,surfhiggs,LT-universal,geomhiggs,
vortexdef}. 
The main steps will be outlined below. 

{\bf Step 1}. 
Since $X$ is a smooth projective curve over $\IC$, it has a 
complex manifold structure, which will be denoted by ${\hat X}$. 
Let $M_1^{\sf an}$, $M_2^{\sf an}$ 
 denote the complex holomorphic 
line bundles on ${\hat X}$ corresponding to the invertible sheaves $M_1, M_2$
on $X$. 
 Let ${\hat M}_1$, ${\hat M}_2$ denote the underlying $C^\infty$-vector bundles 
of the holomorphic line bundles $M_1^{\sf an}$, $M_2^{\sf an}$ and let $\dbar_1$, 
$\dbar_2$ denote the corresponding Dolbeault operators. 
Set ${\hat M}={\hat M}_1\otimes_\IC {\hat M}_2$ and note that there is a 
fixed isomorphism ${\hat M}^\vee{\buildrel \sim \over \longto}  \Lambda^{1,0}_{\hat X}$. 
Let ${\dbar_0}$ denote the canonical Dolbeault operator acting on sections 
of the trivial complex line bundle on ${\hat X}$.

A $V$-framed holomorphic ADHM bundle on ${\hat X}$ 
is defined by the data ${\hat \CE}=({\hat E}, {\dbar}, V, {\hat \Phi}_{1,2}, {\hat \phi}, {\hat \psi})$
where ${\hat E}$ is a $C^\infty$ complex vector bundle on $X$, 
\[
{\dbar}: C^\infty({\hat E}) \to C^\infty({\hat E}\otimes_\IC \Lambda^{0,1}_{\hat X})
\]
 is a semiconnection (or Dolbeault operator) on ${\hat E}$ as defined in 
 \cite[Def. 9.1]{genDTI}, $V$ is a complex vector space of dimension $v\in \{0,1\}$ 
 and 
 \be\label{eq:holADHMa}
 \begin{aligned}
 {\hat \Phi}_i \in C^{\infty}(Hom({\hat E}\otimes_\IC {\hat M}_i , &{\hat E})), \qquad 
 {\hat \phi} \in C^{\infty}(Hom({\hat E}\otimes_\IC {\hat M}, \Lambda^{0,0}_{\hat X}\otimes_\IC V)) \\ & 
 {\hat \psi} \in C^{\infty}({\hat E}\otimes_\IC V^\vee) \end{aligned}
 \ee
 are $C^{\infty}$-morphisms of complex bundles on ${\hat X}$ satisfying the 
 ADHM relation. In addition note that the Dolbeault operators $\dbar$, 
 $\dbar_0$, $\dbar_1$, $\dbar_2$ determine similar differential 
 operator on each space of section in \eqref{eq:holADHMa}, and each morphism
 $({\hat \Phi}_{1,2}, {\hat \phi}, {\hat \psi})$ is required to lie in the 
 kernel of the corresponding Dolbeault operator. Obviously, ${\hat \phi}$ and 
 ${\hat \psi}$ are identically zero if $v=0$. Note that using the fixed isomorphism 
 ${\hat M}^\vee {\buildrel \sim \over \longto} \Lambda^{1,0}_{\hat X}$, 
 the morphism ${\hat \phi}$ is identified with a section in $C^\infty(Hom({\hat E}, 
 \Lambda^{1,0}_{\hat X}\otimes_\IC V))$. Such an identification will be implicit 
 from now on. 
 
 For future reference a $V$-framed $C^\infty$ ADHM bundle is defined 
 by data $({\hat E}, {\hat \Phi}_{1,2}, {\hat \phi}, {\hat \psi})$ as above, 
 except that no holomorphy condition is imposed. The definition of 
 isomorphisms of $C^\infty$ of 
 V-framed ADHM bundles is obvious. 
 
 Given any $C^\infty$ complex vector bundle ${\hat E}$, consider 
 the following infinite dimensional affine space
 \be\label{eq:fieldspace}
\begin{aligned} {\mathscr A}_{ADHM}=\
& {\mathscr A}  \times 
C^{\infty}({Hom}({\hat E}\otimes_\IC ({\hat M}_1\oplus {\hat M}_2), {\hat E})) 
\times \\
& \ C^{\infty}({Hom}({\hat E},  \Lambda^{1,0}_{\hat X}\otimes_\IC V)) 
\times C^{\infty}({\hat E}\otimes_\IC V^\vee)\\
\end{aligned}
\ee 
  where ${\mathscr A}$ is the affine space of semiconnections on ${\hat E}$. 
 
 The group of gauge transformations acting ${\mathscr{A}}_{ADHM}$
 is the infinite dimensional Lie group 
 ${\mathscr G}=C^\infty({\rm{Aut}}({\hat E}))\times (\IC^\times)^{\times v}$ 
 where ${Aut}({\hat E})\subset {End}({\hat E})$
is the subbundle of invertible endomorphisms of ${\hat E}$. 
Note that the second factor $(\IC^\times)^{\times v}$ represents the stabilizer 
of the canonical Dolbeault operator $\dbar_0$ on the trivial line 
bundle $\Lambda^{0,0}_{\hat X} \times_C V$. The later has to be kept fixed in this construction 
since our goal is to 
obtain a local presentation of moduli stacks of objects in the abelian 
category $\CC_X$ rather than the larger abelian category ${\mathcal Q}_X$
of ADHM quiver sheaves on $X$. 

The data $(\dbar, {\hat \Phi}_{1,2}, {\hat \phi}, {\hat \psi})$ will be called simple 
if its stabilizer in ${\mathscr G}$ is isomorphic to the canonical $\IC^\times$ 
subgroup. The subspace of simple data will be denoted by ${\mathscr A}^{\sf si}_{ADHM}$. 

Suppose ${\hat \CE}=({\hat E}, {\hat \Phi}_{1,2}, {\hat \phi}, {\hat \psi}) $ is a $V$-framed holomorphic ADHM bundle on ${\hat X}$ and let 
$(\dbar + A,\Psi_{1,2}, \rho,\eta) \in {\mathscr A}_{ADHM}$,
where $A\in C^{\infty}(End({\hat E})\otimes_\IC \Lambda^{0,1}_{\hat X})$. 
Note that there are natural cup-products 
\be\label{eq:cupprod}
\begin{aligned} 
C^\infty({End}({\hat E})\otimes_\IC \Lambda^{0,1}_{\hat X})\otimes 
C^{\infty}({Hom}({\hat E}\otimes_\IC {\hat M}_i, {\hat E})) & 
\to C^{\infty}({Hom}({\hat E}\otimes_\IC {\hat M}_i, {\hat E})\otimes_\IC \Lambda^{0,1}_{\hat X}) \\
C^\infty({End}({\hat E})\otimes_\IC \Lambda^{0,1}_{\hat X})\otimes 
C^{\infty}({\hat E}\otimes_\IC V^\vee) &  \to C^{\infty}({\hat E}\otimes_\IC V^\vee
\otimes_\IC \Lambda^{0,1}_{\hat X}) \\
C^{\infty}({Hom}({\hat E}, (\Lambda^{1,0}_{\hat X})\otimes_\IC V)) \otimes 
C^\infty({End}({\hat E})\otimes_\IC \Lambda^{0,1}_{\hat X})&\to 
C^{\infty}({Hom}({\hat E}, (\Lambda^{1,1}_{\hat X})\otimes_\IC V))\\
\end{aligned}
\ee
\[\begin{aligned}
C^{\infty}({Hom}({\hat E}\otimes_\IC {\hat M}_i, {\hat E}))
\otimes C^{\infty}({Hom}({\hat E}\otimes_\IC {\hat M}_{3-i}, {\hat E}))
& \to C^{\infty}({End}({\hat E})\otimes_\IC \Lambda^{1,0}_{\hat X}),\\
C^{\infty}({\hat E}\otimes_\IC V^\vee)\otimes 
C^{\infty}({Hom}({\hat E}, (\Lambda^{1,0}_{\hat X})\otimes_\IC V)) 
 &  \to C^{\infty}({End}({\hat E})^{\oplus v}\otimes_\IC \Lambda^{1,0}_{\hat X})\\
 \end{aligned}
\]
where $i=1,2$. Then the data 
\[
({\hat E}, {\overline \partial}+A, {\hat \Phi}_{1}+\Psi_1, {\hat \Phi}_2+\Psi_2, 
{\hat \phi}+\rho, {\hat \psi}+\eta)
\]
defines a $V$-framed  holomorphic  ADHM sheaf structure on 
${\hat E}$ if the following conditions
are satisfied 
\be\label{eq:crtpt}
\begin{aligned} 
{\overline \partial}\Psi_{1} + [A, {\hat \Phi}_1+\Psi_{1}] & =0 \\
{\overline \partial}\Psi_{2} + [A, {\hat \Phi}_2+\Psi_{2}] & =0 \\
{\overline \partial}\rho - (\phi +\rho) A & =0 \\
{\overline \partial}\eta + A (\psi +  \eta) & =0 \\
[{\hat \Phi}_1,\Psi_2] + [\Psi_1, {\hat \Phi}_2] + [\Psi_1,\Psi_2]+\eta\phi +\psi\rho +\eta\rho &=0\\
\end{aligned}
\ee
where all products are the natural cup-products \eqref{eq:cupprod} and 
the commutators are commutators of cup-products. 

Finally, note that 
there is an obvious one-to-one correspondence between $V$-framed holomorphic 
ADHM bundles ${\hat \CE}$ on the complex manifold ${\hat X}$ and locally free objects $\CE$ 
of $\CC_X$ with $v\in \{0,1\}$. 
Moreover, for any $V$-framed holomorphic ADHM bundle ${\hat \CE}$ on ${\hat X}$ 
there is a three term complex $\CC({\hat \CE})$  of $C^\infty$ complex vector
 bundles on ${\hat X}$  constructed by analogy with the deformation complex 
$\CC_{\mathrm{def}}(\CE)$. Given the Dolbeault operator $\dbar$, 
$\CC({\hat \CE})$ admits a canonical 
Dolbeault resolution, which yields in turn a hypercohomology 
double complex. Let  $\CC_{Db}({\hat \CE})$ 
be the total complex of the resulting double complex, and let $C_{Db}({\hat \CE})$ 
denote the complex obtained by taking global $C^\infty$ sections of the terms 
of $\CC_{Db}({\hat \CE})$. Let 
$\IH^k_{Db}({\hat X}, \CC({\hat \CE}))$, $k\in \IZ_{\geq 0}$, 
be the cohomology groups of $C_{Db}({\hat \CE})$. Then there are isomorphisms 
of complex vector spaces 
\[
\IH^k_{Db}({\hat X}, \CC({\hat \CE}))\simeq \IH^{k}(X,\CC(\CE))
\]
for all $k\in \IZ_{\geq 0}$. 

Choosing hermitian structures on ${\hat X},{\hat E}$, the complex $\CC_{Db}({\hat \CE})$ 
will be elliptic, and 
the hypercohomology groups $\IH^k_{Db}({\hat X}, \CC({\hat \CE}))$ are identified 
with spaces of harmonic bundle valued differential forms using Hodge theoretic methods. 
Analogous constructions have been carried in \cite[Sect. 3]{geomhiggs} for Higgs bundles, 
respectively
\cite[Sect. 5]{vortexdef} for triples, hence details will be omitted. 

We can also construct complex Banach manifolds by taking appropriate Sobolev 
completions of the spaces ${\mathscr A}_{ADHM}$, ${\mathscr A}_{ADHM}^{\sf si}$, 
respectively a complex Banach group by taking a Sobolev completion of 
${\mathscr G}$ as in \cite[Sect. 1]{M-Kuranishi}, \cite[Sect. 9.1]{genDTI}.

Given a $C^\infty$ complex vector bundle ${\hat E}$ on ${\hat X}$, 
families of $V$-framed holomorphic ADHM structures on ${\hat E}$ are defined 
by analogy with \cite[Def. 1.5]{M-Kuranishi}, \cite[Def. 9.2]{genDTII}, or 
\cite[Def. 2.3]{KL-moduli}. Then the existence of a versal deformation family 
extending a given $V$-framed holomorphic ADHM bundle ${\hat \CE}$ 
follows from \cite[Thm. 1.1]{locmoduli} by an argument analogous to 
\cite[Thm 2.4]{KL-moduli}. Alternatively, one can check that the proof 
of \cite[Thm 1.]{M-Kuranishi} carries over to the present situation with appropriate 
modifications.

{\bf Step 2.} 
Next note that $X$ also has a structure  of complex analytic space $X^{\sf an}$, and 
one can obviously construct a category of analytic ADHM quiver sheaves 
by analogy with the abelian category $\CC_X$. This category will be denoted 
by $\CC_{X^{\sf an}}$. 
Then given a locally free object $\CE^{\sf an}$ 
of $\CC_{X^{\sf an}}$ with $v\in \{0,1\}$, 
one has to establish the existence 
of a versal deformation family of analytic objects of $\CC_{X^{\sf an}}$ which 
extends $\CE^{\sf an}$. The analogous result for complex analytic vector bundles 
has been proven in \cite{FK,ST-def}.  As observed for example 
in \cite{KL-moduli} the extension to families of decorated 
analytic bundles follows from the complex analytic version of the standard representability 
for Hom functors presented for example in \cite[Thm 5.8]{FGAexplained}. 
The complex analytic version of this result follows from  \cite{flenner}. 

Now let ${\hat \CE}$ be a $V$-framed holomorphic ADHM bundle on ${\hat X}$
and let $T$ denote the base of the versal deformation family of $V$-framed 
holomorphic ADHM bundles extending ${\hat \CE}$; $T$ is a finite 
dimensional complex analytic space. 
Let $\pi_{X^{\sf an}}: X^{\sf an}\times T \to T$ denote the canonical projection.   
Let ${\CE}^{\sf an}$
the locally free object of $\CC_{X^{\sf an}}$ corresponding to 
${\hat \CE}$. 
Then by analogy with \cite[Prop 2.3]{M-Kuranishi}, \cite[Prop. 9.5]{genDTI}, 
or \cite{KL-moduli}[Thm. 2.5] there exists a versal  deformation family $\CE^{\sf an}_T$ 
with base $T$ extending  ${\CE}^{\sf an}$ and an isomorphism 
$\CE^{\sf an}_T{\buildrel \sim \over \longto} \pi_{X^{\sf an}}^* {\hat \CE}$ 
of $C^\infty$ $V$-framed ADHM sheaves which induces the versal deformation family 
of ${\hat \CE}$. Moreover, $\CE^{\sf an}_T$ is universal if ${\hat \CE}$ is simple. 

{\bf Step 3.} Next let $\CE$ be a locally free object of $\CC_X$ and 
let $\CE^{\sf an}$ be the corresponding complex analytic object. 
Using the standard representability result \cite[Thm 5.8]{FGAexplained}
for Hom functors, and the existence of an algebraic versal deformation family 
extending $\CE$ proven in \cite[Prop. 9.8]{genDTI}, it follows that 
there exists an algebraic versal deformation family of objects of $\CC_X$ 
extending $\CE$. Moreover, \cite[Prop. 9.9]{genDTI} proves that 
the algebraic and the analytic versal deformation families 
associated to a given holomorphic vector bundle are locally isomorphic  with respect to the 
complex analytic topology. The extension of this result to ADHM sheaves 
is straightforward.  

{\bf Step 4.} In order to conclude the proof of theorems (\ref{theoremA}, \ref{theoremB}) 
it suffices to prove the existence of a holomorphic functional on the Sobolev completion 
of the space ${\mathscr A}_{ADHM}$, for any $V$-framed holomorphic 
ADHM sheaf ${\hat \CE}$ on ${\hat X}$,  
with the same properties as the holomorphic 
Chern-Simons functional employed in \cite[Sect. 9.5]{genDTI}. 

As observed in step 1 above, 
the Dolbeault operators $\dbar$, $\dbar_0$, $\dbar_1$, $\dbar_2$ 
determine similar Dolbeault operators on all spaces of sections in the right hand side of 
\eqref{eq:fieldspace}. In order to keep the notation short all the resulting operators 
will be denoted by the same symbol $\dbar$, the distinction being clear once the 
argument of the operator is specified. 

Then the  holomorphic Chern-Simons functional for $V$-framed ADHM sheaves is defined by 
\be\label{eq:holCSact}
\begin{aligned}
CS: (A,\Psi_{1,2}, \rho,\eta) \to 
\int_{X} \mathrm{Tr}( & \Psi_2 {\overline \partial} \Psi_1 + \rho \dbar\eta + 
A[{\hat \Phi}_1,\Psi_2] + 
A[\Psi_1,{\hat \Phi}_2] + A[\Psi_1,\Psi_2] \\
& + A\psi\rho + A \eta \phi + A \eta\rho)\\
\end{aligned}
\ee
where all products are the natural cup-products \eqref{eq:cupprod} and 
the commutators are commutators of cup-products. Note that given the cup-products 
\eqref{eq:cupprod} it is straightforward to check that the integrand in the right hand side of 
equation \eqref{eq:holCSact} is a section of $\Lambda^{1,1}_X$. It is also straightforward 
to check that the functional \eqref{eq:holCSact} is gauge invariant, and 
the critical points of the functional \eqref{eq:holCSact} are determined by the
equations \eqref{eq:crtpt}. 

\hfill $\Box$

\begin{lemm}\label{movingeulerlemma}
Let $\CE_1,\CE_2$ be two locally free ADHM sheaves on $X$ of numerical types  
$(r_1,e_1,1)$, $(r_2,e_2,0)$. Let 
\[
\begin{aligned}
\chi(\CE_1,\CE_2) = &\ \mathrm{dim}\mathrm{Ext}^{0}(\CE_1,\CE_2)-
\mathrm{dim}\mathrm{Ext}^{1}(\CE_1,\CE_2) \\
& - \mathrm{dim}\mathrm{Ext}^{0}(\CE_2,\CE_1) + 
\mathrm{dim}\mathrm{Ext}^{1}(\CE_2,\CE_1).\
\end{aligned} 
\]
Then $\chi(\CE_1,\CE_2)= e_2-r_2(g-1)$ depends only on the 
numerical types of $\CE_1,\CE_2$. 
\end{lemm}

{\it Proof.} Recall that corollary (\ref{Cextcoro}) proves that the extension groups 
$\mathrm{Ext}^1(\CE_1,\CE_2)$ are the hypercohomology groups of 
the three term locally free complex $\CC(\CE_1,\CE_2)$ written in 
equation \eqref{eq:hypercohA}. 
Then lemma (\ref{movingeulerlemma}) follows by a straightforward 
application of  the Riemann-Roch theorem. 

\hfill $\Box$

Next let $\nu$ denote Behrend's constructible function for the Artin stack 
$\obj(\CX)_{\leq 1}$ constructed in  \cite[Prop. 4.4]{genDTI}.
Then the following theorem holds by analogy with \cite[Thm. 5.9]{genDTI}. 
\begin{theo}\label{nuidentity} 
Let $\CE_1,\CE_2$ be locally free objects of $\CC_X$ with $v(\CE_1)+v(\CE_2)\leq 1$.
Then the following identities hold 
\be\label{eq:nuidA}
\begin{aligned}
\nu([\CE_1\oplus \CE_2]) = (-1)^{\chi(\CE_1,\CE_2)} \nu([\CE_1]) \nu([\CE_2)])
\end{aligned}
\ee
 \be\label{eq:nuidAB}
\begin{aligned}
& \mathop{\int_{[\CE]\in \IP(\mathrm{Ext}^1(\CE_2,\CE_1))}}_{0\to \CE_1\to \CE\to \CE_2\to 0}  \nu([\CE]) d\chi - 
\mathop{\int_{[\CE]\in \IP(\mathrm{Ext}^1(\CE_1,\CE_2))}}_{0\to \CE_2\to \CE\to \CE_1\to 0} \nu([\CE]) d\chi = \\
& \qquad \mathrm{dim}\mathrm{Ext}^{1}(\CE_2,\CE_1)
- \mathrm{dim}\mathrm{Ext}^{1}(\CE_1,\CE_2)\\
\end{aligned}
\ee
\end{theo}

 \appendix 
\section{Higgs sheaves}\label{higgsapp} 
For completeness we summarize the main results on moduli 
of Higgs sheaves used 
in this paper. 
Let $X$ be a smooth projective curve over a field $K$ 
over $\IC$. Let $(M_1,M_2)$ be fixed line bundles 
on $X$. 

\begin{defi}\label{Higgs-sheaves}
$(i)$ A Higgs sheaf of type $(r,e)$ on $X$ with coefficient sheaves
$(M_1,M_2)$
is a collection $\CE=(E,\Phi_1,\Phi_2)$ where $E$ is a coherent 
sheaf of type $(r,e)$ on $X$ and 
\[
\Phi_{i}:E\otimes_X M_i \to E
\]
are morphisms of $\CO_X$-modules satisfying 
\be\label{eq:higgsrel}
\Phi_1\circ (\Phi_2\otimes 1_{M_1}) - \Phi_2\circ (\Phi_1\otimes 
1_{M_2}) =0.
\ee

$(ii)$ A morphism of Higgs sheaves $\CE,\CE'$ is a morphism 
$\xi:E\to E'$ of coherent sheaves on $X$ satisfying the obvious 
compatibility conditions with the data $\Phi_{1,2}$, $\Phi_{1,2}'$. 

$(iii)$ A Higgs sheaf $\CE=(E,\Phi_1,\Phi_2)$
of type $(r,e)$ on $X$  is called (semi)stable if any nontrivial 
$\Phi$-invariant proper 
saturated subsheaf $0\subset E'\subset E$ satisfies 
\be\label{ssHiggs}
\mu(E') \ (\leq) \ \mu(E)
\ee
\end{defi}

The following theorem summarizes the properties of moduli of Higgs sheaves 
following \cite{semistpairs,simpsonI,projective}. 
\begin{theo}\label{higgsmoduli}
Suppose $X$ is a smooth projective curve over $\IC$ and let 
$M_1,M_2$ be fixed line bundles on $X$ as in the main text. 
Let $\CX=(X,M_1,M_2)$. Then 

$(i)$ For $(r,e)\in \IZ_{\geq 1}\times \IZ$ 
there is a quasi-projective coarse 
moduli scheme $Higgs^{ss}(\CX,r,e)$ over $\IC$
parameterizing 
$S$-equivalence classes of Higgs sheaves of type $(r,e)$ on $X$ 
with coefficient sheaves $(M_1,M_2)$. The scheme $Higgs^{ss}(\CX,r,e)$ contains 
an open subscheme $Higgs^{s}(\CX,r,e)$ which parameterizes isomorphism classes 
of stable Higgs sheaves. 

$(ii)$ There is a proper Hitchin morphism 
$h: Higgs^{ss}(\CX,r,e) \to \IH$ where 
 $\IH = \oplus_{n=0}^r H^0(X, \mathrm{Symm}^n(M_1^{-1}\oplus M_2^{-1}) )$
 mapping a polystable Higgs sheaf of type $(r,e)$ to its characteristic polynomial. 
\end{theo}

\bibliography{newref.bib}
 \bibliographystyle{abbrv}
\end{document}